\documentclass[10pt]{amsart}
\usepackage{amsmath,amstext,amsfonts,%
enumitem,marvosym,%
extpfeil,graphicx,mathrsfs,%
mathtools,setspace,tikz,xparse}
\usepackage[alphabetic]{amsrefs}
\usepackage[unicode]{hyperref}
\usepackage[utf8]{inputenc}
\usepackage[T1]{fontenc}
\usepackage[normalem]{ulem}
\usepackage[cmtip,all]{xy}
\usepackage[light]{sourcesanspro}
\usepackage{lmodern}
\hypersetup{colorlinks=false,linkbordercolor={0.93 0.71 0.13},citebordercolor={1 0 .5},urlbordercolor={0.27 0.39 0.68}}%
\usetikzlibrary{decorations.text,shapes.geometric,arrows,positioning,fit,calc,matrix}

\usepackage{filecontents}

\definecolor{MyColor1}{HTML}{fff3b0}
\definecolor{MyColor2}{HTML}{e09f3e}
\definecolor{MyColor3}{HTML}{335c67}
\definecolor{MyColor4}{HTML}{393e41}

\DeclareMathOperator{\Ker}{ker}

\DeclareMathOperator{\Img}{im}

\DeclareMathOperator{\ind}{ind}

\newcommand{\m}[1]{\(\smash{#1}\)}

\newcommand{\bb}[1]{\mathbb{#1}}
\newcommand{\QQ}{\bb{Q}}
\newcommand{\FF}{\bb{F}}
\newcommand{\ZZ}{\bb{Z}}

\newcommand{\FFp}{\overline{\bb{F}}_p}
\newcommand{\QQp}{\smash{\overline{\bb{Q}}_p}}

\newcommand{\ZZp}{\smash{\overline{\bb{Z}}_p}}

\newcommand{\defeq}{\vcentcolon=}
\newcommand{\SCR}{\mathscr}
\newcommand{\IE}{i.e.\ }

\newcommand{\HECKET}{\SCR T}

\newcommand{\MU}{\mu}

\newcommand{\SIGMA}{\Sigma}

\newcommand{\LAMBDA}{\lambda}

\newcommand{\LESS}[1]{\smash{\uline{#1}}}
\newcommand{\MORE}[1]{\smash{\overline{#1}}}
\newcommand{\BIGO}[1]{\smash{\mathsf{O}\!\left({#1}\right)}}

\newcommand{\TEMPORARYA}{\xi}

\newcommand{\TEMPORARYD}{\vartheta}%

\newcommand{\TEMPORARYI}{\nu}
\newcommand{\TEMPORARYM}{\mathrm{M}}

\newcommand{\PADICVALUATION}{v_p}

\newcommand{\TEMPORARYF}{L}
\newcommand{\TEMPORARYDISC}{D}

\newcommand{\sfirst}{s_{1}}
\newcommand{\ssecond}{s_{2}}
\newcommand{\TILDE}{\widetilde}

\newcommand{\GEQSLANT}{\geqslant}
\newcommand{\LEQSLANT}{\leqslant}

\newcommand{\maxideal}{\mathfrak m}

\newcommand{\NEWLINEA}{\\[2pt]}

\newcommand{\SUBMODULE}{\mathrm{sub}}
\newcommand{\QUOTIENTI}{\mathrm{quot}}

\newcommand{\FILT}{\widehat}
\newcommand{\FORMX}{X}
\newcommand{\FORMY}{Y}
\newcommand{\DELTAA}{\Delta}
\newcommand{\TEMPORARYJ}{\eta}

\newcommand{\sqbrackets}[4]{{#1}\bullet_{#2,#3}{#4}}

\newtheorem{theorem}{Theorem}
\newtheorem{lemma}[theorem]{Lemma}
\newtheorem{proposition}[theorem]{Proposition}
\newtheorem{corollary}[theorem]{Corollary}
\newtheorem{conjecture}{Conjecture}

\newenvironment{Smatrix}%
    {\tiny (\begin{smallmatrix}}%
    {\end{smallmatrix})}
\newenvironment{Lmatrix}%
    {\footnotesize (\begin{smallmatrix}}%
    {\end{smallmatrix})}
\newcommand{\BLACKSQUARE}{\hspace{\stretch{1}}\m{\mathbin{\vcenter{%
    \hbox{\rule{.9ex}{.9ex}}}}}}
\renewenvironment{proof}[1]%
    {\emph{{#1}}}%
    {\BLACKSQUARE}

\makeatletter
\newcommand{\xhookdoubleheadrightarrow}[2][]{%
  \lhook\joinrel
  \ext@arrow 0359\rightarrowfill@ {#1}{#2}%
  \mathrel{\mspace{-15mu}}\rightarrow}
\renewcommand{\xtwoheadrightarrow}[2][]{%
  \ext@arrow 0359\rightarrowfill@ {#1}{#2}%
  \mathrel{\mspace{-15mu}}\rightarrow}
\def\@settitle{\begin{center}%
  \baselineskip14\p@\relax\normalfont\LARGE
    \uppercasenonmath\@title
  \@title\end{center}%
}
\def\do#1{\@ifdefinable#1{\newdimen#1}}
\do\@TempWidthOne
\do\@TempWidthTwo
\newcommand*\WidthInMathUnitsOf[2]{%
    \settowidth\@TempWidthOne{$\m@th #1#2$}%
    \settowidth\@TempWidthTwo{$\m@th #1\mspace{1mu}$}%
    \strip@pt \dimexpr \@TempWidthOne*\p@/\@TempWidthTwo \relax
}
\newcommand*\SetToWidthInMathUnits[3]{%
    \settowidth\@TempWidthOne{$\m@th #2#3$}%
    \settowidth\@TempWidthTwo{$\m@th #2\mspace{1mu}$}%
    \edef #1{%
        \strip@pt \dimexpr \@TempWidthOne*\p@/\@TempWidthTwo \relax
    }%
}
\newcommand*\ExtractHexDigitOfFamilyOfSymbol[1]{%
    \expandafter\@ExtractHexDigitOfFamily \meaning #1%
}
\@ifdefinable\@ExtractHexDigitOfFamily{
    \def\@ExtractHexDigitOfFamily#1"#2#3#4#5{#3}
}
\makeatother

\newcommand*{\ReportValueOfOneEm}[1]{%
    &\text{in #1 size,}%
        &\quad\texttt{1em}%
            &=\texttt{\the \fontdimen 6 \csname #1font\endcsname 2}%
        &&\qquad\text{(from \texttt{\fontname \csname #1font\endcsname 2})}%
}
\newcommand*{\ReportSizeOfFontOfSymbol}[2]{%
    &\texttt{\pdffontsize\csname #2font\endcsname
                "\ExtractHexDigitOfFamilyOfSymbol{#1}}%
        &&\quad\text{in #2 size}%
        &\qquad(\texttt{\string\fontname}%
            &=\texttt{\fontname\csname #2font\endcsname
                "\ExtractHexDigitOfFamilyOfSymbol{#1}})%
}

\title%
[Reductions of crystalline representations, III]%
{On the reductions of certain
two-dimensional crystalline representations, III}
\author{Bodan Arsovski}
\date{February 2021}

\begin{document}
\maketitle

\renewcommand{\HECKET}{T}
\newcommand{\RRRRRR}{\SCR R^{s,\TEMPORARYI}}
\newcommand{\VERIFYIF}[1]{[{#1}]}
\newcommand{\DARKCIRC}{\bullet}
\newcommand{\TEMPORARYAAAA}{w}
\newcommand{\UUU}{U}
\newcommand{\subm}{sub}

\newcommand{\MATRIXNO}{\SCR M_1}
\newcommand{\MATRIXYES}{\SCR M_2}
\renewcommand{\TEMPORARYF}{m}
\renewcommand{\TEMPORARYDISC}{\SCR N}
\newcommand{\LARGESTSLOPE}{100}
\newcommand{\LARGESTPRIME}{139522597633}
\newcommand{\ddd}{\mathrm{d}}
\newcommand{\polyr}{z}
\newcommand{\polyH}{\Phi}
\newcommand{\polys}{\psi}
\newcommand{\polyconst}{C}
\newcommand{\veeC}{\check{C}}
\newcommand{\vveeC}{\hat{C}}
\newcommand{\annd}{\text{ \& }}
\newcommand{\Cnegativeone}{\frac{(-1)^{\beta}(s-\alpha)(\alpha-\beta+1)}{\beta^2(2\alpha-s+1) \binom{\alpha}{\beta}} \binom{\alpha}{s-\alpha}}
\newcommand{\rowR}{\mathbf{r}}
\newcommand{\colA}{\mathbf{c}}
\newcommand{\abcNL}{\\[3pt]}
\newcommand{\Cjnewline}{\\[5pt]}
\newcommand{\tempheightaer}{\rule{0pt}{12pt}}
\newcounter{equ}
\newcommand{\jj}{t}
\newcommand{\llll}{\gamma}
\newcommand{\zs}{u}
\newcommand{\za}{v}
\newcommand{\zb}{w}
\newcommand{\zi}{t}
\newcommand{\varu}{e}
\newcommand{\boxcol}[1]{\tikz{\path[draw=black,fill={#1}] (0,0) rectangle (0.2cm,0.2cm);}}

\begin{abstract}
A conjecture of Breuil, Buzzard, and Emerton says that the slopes %
 of certain  %
 reducible \m{p}-adic Galois representations must be integers. In~\cite{Ars1} we showed this conjecture for %
  representations that lie over certain ``non-subtle'' components of weight space. %
   This article is a continuation of~\cite{Ars1} in which we also show the conjecture for the ``subtle'' components for slopes less than \m{\frac{p-1}{2}}.
\end{abstract}

\section{Introduction and results}

\subsection{Background} Let \m{p} be an odd prime number %
 and \m{k\GEQSLANT2} be an integer, and let \m{a} be an element of  \m{\ZZp} such that \m{\PADICVALUATION(a)>0}. Let us denote \m{\TEMPORARYI = \lfloor \PADICVALUATION(a) \rfloor + 1 \in \ZZ_{>0}}. With this data one can associate a certain two-dimensional crystalline \m{p}-adic representation \m{V_{k,a}} with Hodge--Tate weights \m{(0,k-1)}. We give the definition of this representation in section~2 of~\cite{Ars1}, and we define \m{\overline{V}_{k,a}} as the semi-simplification of the reduction modulo the maximal ideal \m{\maxideal} of \m{\ZZp} of a Galois stable \m{\ZZp}-lattice in \m{{V_{k,a}}} (with the resulting representation being independent of the choice of lattice). The question of computing \m{\overline{V}_{k,a}} has been studied extensively, and we refer to the introduction of~\cite{Ars1} for a brief exposition of it. Partial results have been obtained by %
  Fontaine, Edixhoven, Breuil, Berger, Li, Zhu, Buzzard, Gee, Bhattacharya, Ganguli, Ghate, et al (see~\cite{Ber10}, \cite{Bre03a}, \cite{Bre03b}, \cite{Edi92}, \cite{BLZ04}, \cite{BG15}, \cite{BG09}, \cite{BG13}, \cite{GG15}). 
A conjecture of Breuil, Buzzard, and Emerton says the following.
\begin{conjecture}\label{MAINCONJ}
{If \m{k} is even and \m{\PADICVALUATION(a) \not\in \ZZ} then \m{\overline{V}_{k,a}} is irreducible.}
\end{conjecture}

The main result of~\cite{Ars1} is that this conjecture is true over
certain ``non-subtle'' components of weight space.
We say that a weight \m{k} belongs to a  ``non-subtle'' component of weight space if and only if
\[\textstyle \textstyle k \not\equiv 3,4,\ldots,2\TEMPORARYI,2\TEMPORARYI+1 \bmod p-1.\]
Thus there are \m{\max\{\frac{p-1}{2} - \TEMPORARYI+1,0\}} many ``non-subtle'' components of weight space. 
This article is a continuation of~\cite{Ars1} in which we also show the conjecture for the ``subtle'' components for slopes less than \m{\frac{p-1}{2}}.
The main result we show is the following theorem.

\begin{theorem}\label{MAINTHM1}
Conjecture~\ref{MAINCONJ} is true when the slope is less than \m{\frac{p-1}{2}}.
\end{theorem}

    \section{Computing \texorpdfstring{\m{\overline V_{k,a}}}{V} by computing \texorpdfstring{\m{\overline\Theta_{k,a}}}{\textTheta}}
    \label{secComp}
    
    From now on we assume the notation from sections 2, 3, 4, and 6 of~\cite{Ars1}. Moreover, we assume that \m{k>p^{100}} as in section~5 of~\cite{Ars1}.
          Theorem~2 in~\cite{Ars1} implies that our main theorems can be rewritten in the following equivalent forms. Recall that we assume \m{p>2} throughout.

\begin{theorem}\label{MAINTHM1a}
If \m{k\in 2\ZZ} and \m{\PADICVALUATION(a) \in (0,\frac{p-1}{2}) \backslash \ZZ}  then \m{\overline{\Theta}_{k,a}} is irreducible.
\end{theorem}
    
Thus our task is to prove theore~\ref{MAINTHM1a}.

\section{Combinatorics}
\label{secComb}

Throughout the proof we will refer to the combinatorial results in section 8 of~\cite{Ars1}. For convenience, we reproduce the statements here in the form we will use.

\begin{lemma}\label{part1X5}
Suppose throughout this lemma that 
\begin{align*}
    n,t,y \in \ZZ, \qquad
    b,d,k,l,w  \in \ZZ_{\GEQSLANT 0}, \qquad
    m,u,v  \in \ZZ_{\GEQSLANT 1}.
\end{align*}
\begin{enumerate}
    \item If \m{u \equiv v \bmod (p-1)p^{m-1}} then
\begin{align}
\label{ceqn1}\tag{\emph{c-a}}
\textstyle \TEMPORARYM_{u,n} \equiv \TEMPORARYM_{v,n} \bmod p^m.
\end{align}
\item Suppose that \m{u = t_u(p-1)+s_u} with \m{s_u = \MORE{u}}, so that \m{s_u \in \{1,\ldots,p-1\}} and \m{t_u \in \ZZ_{\GEQSLANT 0}}. Then \begin{align}
\label{ceqn2}\tag{\emph{c-b}}
\textstyle \TEMPORARYM_u = 1 + \VERIFYIF{u\equiv_{p-1} 0} + \frac {t_u}{s_u}p + \BIGO{t_up^2}.
\end{align}
\item  If \m{n\LEQSLANT 0} then
\begin{align}
\label{ceqn3}\tag{\emph{c-c}}
\textstyle \TEMPORARYM_{u,n} = %
\sum_{i=0}^{-n} (-1)^i \binom{-n}{i} \TEMPORARYM_{u-n-i,0}. %
\end{align} \item If \m{n\GEQSLANT 0} then
\begin{align}
\label{ceqn4}\tag{\emph{c-d}}\textstyle \TEMPORARYM_{u,n} \equiv (1+\VERIFYIF{u\equiv_{p-1}n\equiv_{p-1} 0})\binom{\MORE{u}}{\LESS{n}} %
 \bmod p.
\end{align}
\item If \m{u \GEQSLANT (b+l)d} and \m{l\GEQSLANT w} then
\begin{align}
\label{ceqn5}\tag{\emph{c-e}}\textstyle \sum_j (-1)^{j-b} \binom{l}{j-b}\binom{u-dj}{w} = \VERIFYIF{w=l}d^l.
\end{align}
\item If \m{\FORMX} is a formal variable  then
\begin{align}\label{ceqn8}\tag{\emph{c-f}}
   \textstyle  \binom{\FORMX}{t+l} \binom{t}{w} = \sum_v (-1)^{w-v} \binom{l+w-v-1}{w-v} \binom{\FORMX}{v} \binom{\FORMX-v}{t+l-v}.
\end{align}
Consequently, if \m{b+l\GEQSLANT d+w} then
\begin{align}
\label{ceqn6}\nonumber& \textstyle \sum_i %
 \binom{b-d+l}{i(p-1)+l}\binom{i(p-1)}{w} \\ & \textstyle \null\qquad = \sum_v (-1)^{w-v} \binom{l+w-v-1}{w-v} \binom{b-d+l}{v} \TEMPORARYM_{b-d+l-v,l-v}.\tag{\emph{c-g}}\end{align}
\item We have \begin{align}
\label{ceqn9}\tag{\emph{c-i}}
\textstyle \sum_{j} (-1)^j \binom{y}{j} \binom{y+l-j}{w-j} = (-1)^w \binom{w-l-1}{w}.
\end{align}
\item We have 
\begin{align}
\label{ceqn10}\tag{\emph{c-j}}
\textstyle \sum_{j} \binom{u-1}{j-1} \binom{-l}{j-w} = (-1)^{u-w} \binom{l-w}{u-w}.
\end{align}
\item We have
\begin{align}
\label{ceqn11}\tag{\emph{c-k}}
\textstyle \sum_{j} (-1)^{j} \binom{j}{b} \binom{l}{j-w} = (-1)^{l+w}\binom{w}{l+w-b}.
\end{align}
\end{enumerate}
\end{lemma}

\begin{lemma}\label{part2X5}
Let \m{\alpha \in \ZZ \cap [0,\ldots,\frac{r}{p+1}]} and let \m{\{D_i\}_{i\in \ZZ}} be a family of elements of \m{\ZZ_p} such that \m{D_i=0} for \m{i\not\in[0,\frac{r-\alpha}{p-1}]} and
\m{\textstyle \TEMPORARYD_{w} (D_\DARKCIRC) = 0} for all \m{0 \LEQSLANT w<\alpha}. Then \[\textstyle \sum_i D_i x^{i(p-1)+\alpha}y^{r-i(p-1)-\alpha} = \theta^\alpha h\] for some polynomial \m{h} with integer coefficients.
\end{lemma}

\begin{lemma}\label{part3X5} For \m{\alpha,\lambda,\mu \in \ZZ_{\GEQSLANT 0}} let
\[\textstyle L_\alpha(\lambda,\mu)\] be the 
\m{(\alpha+1)\times(\alpha+1)} matrix with %
 entries
    \[\textstyle L_{l,j} = \sum_{k=0}^{\alpha} \frac{j!}{l!} \left(\frac\mu\lambda\right)^k  \sfirst (l,k) \ssecond (k,j),\]
where \m{\sfirst (l,k)} are the Stirling numbers of the first kind and \m{\ssecond (k,j)} are the Stirling numbers of the second kind. Then
\[\textstyle L_\alpha(\lambda,\mu) \left(\binom{\lambda \FORMX}{0},\ldots,\binom{\lambda \FORMX}{\alpha}\right)^T = \left(\binom{\mu \FORMX}{0},\ldots,\binom{\mu \FORMX}{\alpha}\right)^T.\]
\end{lemma}

\begin{lemma}\label{lemmaComb2} 
For \m{\alpha \in \ZZ_{\GEQSLANT 0}} let \m{B_{\alpha}} be the 
\m{(\alpha+1)\times(\alpha+1)} matrix with %
 entries
    \[\textstyle B_{i,j} = j! \sum_{k,l=0}^{\alpha} \frac{(-1)^{i+l+k}}{l!} \binom{l}{i} (1-p)^{-k} \sfirst (l,k) \ssecond (k,j),\]
where \m{\sfirst (i,j)}  and \m{\ssecond (k,j)} are the Stirling numbers of the first and second kind, respectively. Let \m{\{X_{i,j}\}_{i,j\GEQSLANT 0}} be formal variables.
For \m{\beta \in \ZZ_{\GEQSLANT 0}} such that \m{\alpha\GEQSLANT \beta} let  \[\textstyle S({\alpha,\beta})=\left(S({\alpha,\beta})_{w,j}\right)_{0\LEQSLANT w,j\LEQSLANT\alpha}\] be the \m{(\alpha+1)\times(\alpha+1)} matrix with %
 entries
 \[\textstyle S({\alpha,\beta})_{w,j} = \sum_{i=1}^\beta X_{i,j} \binom{i(p-1)}{w}.\]
 Then \m{B_\alpha S({\alpha,\beta})} is zero outside the rows indexed \m{1,\ldots,\beta} and \[\textstyle (B_\alpha S({\alpha,\beta}))_{i,j} = X_{i,j}\]
for \m{i \in \{1,\ldots,\beta\}}.
\end{lemma}

\begin{lemma}\label{part2X3a}
For \m{u,v,c \in \ZZ} let us define
\[\textstyle F_{u,v,c}(X) = \sum_w (-1)^{w-c} \binom{w}{c} \binom{X}{w}^\partial \binom{X+u-w}{v-w} \in \QQ_p[X].\]
Then
\[\textstyle F_{u,v,c}(X) =  \binom{u}{v-c} \binom{X}{c}^\partial-\binom{u}{v-c}^\partial \binom{X}{c}.\]
\end{lemma}

\begin{lemma}\label{part4X5} Let \m{\FORMX} and \m{\FORMY} denote formal variables, and let
\[\textstyle c_j = (-1)^j \alpha! \left( \frac{\FORMX+j+1}{j+1}\binom{\FORMY}{\alpha-j-1} + \binom{\FORMY}{\alpha-j}\right)\in\QQ[\FORMX,\FORMY]\subset \QQ(\FORMX,\FORMY)\]
be polynomials over \m{\QQ} of degrees \m{\alpha-j}, for \m{1\LEQSLANT j\LEQSLANT \alpha}. Let \[\textstyle M = (M_{w,j})_{0\LEQSLANT w,j\LEQSLANT \alpha}\] be the \m{(\alpha+1)\times(\alpha+1)} matrix over \m{\QQ(\FORMX,\FORMY)} with entries
 \begin{align*}
        M_{w,0} & =\textstyle (-1)^w\frac{(\FORMY-\FORMX)\FORMX_w}{\FORMY_{w+1}}, %
        \\ M_{w,j} & = \textstyle \sum_v (-1)^{w-v} \binom{j+w-v-1}{w-v} \binom{\FORMX+j}{v} \left(\binom{\FORMY+j-v}{j-v}-\binom{\FORMX+j-v}{j-v}\right), %
\end{align*}
for \m{0\LEQSLANT w\LEQSLANT \alpha} and \m{0<j\LEQSLANT\alpha}. Then the first \m{\alpha-1} entries of
\[\textstyle Mc  = M(\FORMY_{\alpha},c_1,\ldots,c_\alpha)^T = (d_0,\ldots,d_\alpha)^T\] are zero, and \m{d_\alpha=\frac{(\FORMY-\FORMX)_{\alpha+1}}{\FORMY-\alpha}}. \end{lemma}

\begin{lemma}\label{lemmaComb4} 
Suppose that \m{s,\alpha,\beta\in\ZZ} are such that
\[\textstyle 1 \LEQSLANT \beta \LEQSLANT \alpha \LEQSLANT \frac{s}{2}-2 \LEQSLANT \frac{p-5}{2}.\]
Let \m{B=B_\alpha} denote the matrix defined in lemma~\ref{lemmaComb2}.
Let \m{M} denote the \m{(\alpha+1)\times(\alpha+1)} matrix with entries in \m{\FF_p} such that if \m{i\in\{1,\ldots,\beta\}} and \m{j\in\{0,\ldots,\alpha\}} then
\[\textstyle M_{i,j} = \binom{\beta}{i} \cdot \left\{\begin{array}{rl}
    \binom{s-\alpha-\beta+i}{i}^{-1} (-1)^{i+1} & \text{if } j = 0, \\
    \binom{s-\alpha-\beta+j}{j-i} & \text{if }j>0,
\end{array}\right.\]
and if \m{i\in\{0,\ldots,\alpha\}\backslash\{1,\ldots,\beta\}} and \m{j\in\{0,\ldots,\alpha\}} then \m{M_{i,j}} is the reduction modulo \m{p} of
 \begin{align*}
       & \textstyle p^{-\VERIFYIF{j=0}}  \sum_{w=0}^\alpha B_{i,w} \sum_v (-1)^{w-v}  \binom{j+w-v-1}{w-v} \binom{s+\beta(p-1)-\alpha+j}{v}^\partial  
       \\
        & \textstyle \null\qquad\qquad\qquad\qquad\qquad\qquad\qquad\qquad\qquad\qquad \cdot\sum_{u=0}^{\beta} \binom{s+\beta(p-1)-\alpha+j-v}{u(p-1)+j-v}
        \\
        & \textstyle \null\qquad - \VERIFYIF{i=0} p^{-\VERIFYIF{j=0}} \textstyle  \binom{s+\beta(p-1)-\alpha+j} {j}^\partial \\
        & \textstyle \null\qquad  -\VERIFYIF{j=0} \sum_{w=0}^\alpha B_{i,w} (-1)^{w}\binom{s+\beta(p-1)-\alpha}{w}  \frac{w!}{(s-\alpha)_{w+1}}.
\end{align*}
Then there is a solution of
\[M (z_0,\ldots,z_\alpha)^T = (1,0,\ldots,0)^T\] such that \m{z_0 \ne 0}. 
\end{lemma}

Now let us prove some additional combinatorial results.

\begin{lemma}\label{lemmaComb5}
Suppose that \m{s,\alpha\in\ZZ_{\GEQSLANT0}} are such that
\[\textstyle  s \in \{2,4,\ldots,p-3\} \text{ and } \frac s2 \LEQSLANT \alpha < s \text{ and } \alpha \LEQSLANT \frac{p-3}{2}.\] For \m{w,j\in\ZZ_{\GEQSLANT 0}} let \m{F_{w,j}(\polyr) \in \FF_p[\polyr]} denote the polynomial
\[\textstyle \sum_v (-1)^{w-v} \binom{j+w-v-1}{w-v}\binom{\polyr-\alpha+j}{v} \binom{\MORE{s-\alpha+j-v}}{\LESS{j-v}} - \binom{\polyr-\alpha+j}{j}\binom{0}{w} - \binom{\polyr-\alpha+j}{s-\alpha}\binom{\polyr-s}{w}.\]
Let  \m{\polyconst_{0}(\polyr),\ldots,\polyconst_\alpha(\polyr) \in \FF_p[\polyr]} denote the polynomials
\begin{align*} 
 \polyconst_j(\polyr) = \left\{\begin{array}{rl} %
 \binom{\alpha}{s-\alpha-1}^{-1}\frac{s-\polyr}{\alpha+1} %
& \text{if } j=0,
\\[5pt] \frac{(-1)^{j+1}}{j+1} \binom{s-\alpha-1}{\alpha-j} (\polyr-\alpha) & \text{if } j \in \{1,\ldots,\alpha\}. \end{array}\right.\end{align*}
Let \m{F_1(\polyr),F_2(\polyr) \in \FF_p[\polyr]} denote the polynomials
 \begin{align*}
     F_1(\polyr) & \textstyle = \sum_{j=0}^\alpha \polyconst_j(\polyr) F_{w,j}(\polyr), \\
     F_2(\polyr) & \textstyle =-\VERIFYIF{w=0} \binom{s-\polyr-1}{\alpha+1}.
 \end{align*}
Note that all of these polynomials depend on \m{s} and \m{\alpha}. Then \m{F_1(z) = F_2(z)}.
\end{lemma}

\begin{proof}{Proof. }
Let us first show that
\begin{align}\label{eqn423}%
\textstyle \polyconst_0(\polyr) \binom{\polyr-\alpha}{s-\alpha} + \sum_{j=1}^\alpha \polyconst_j(\polyr) \binom{\polyr-\alpha+j}{s-\alpha} = \frac{(-1)^{\alpha+1}(\polyr-\alpha)}{s-\alpha}.\end{align}
 Since \[\textstyle \polyconst_0(\polyr) \binom{\polyr-\alpha}{s-\alpha} = \binom{\alpha}{s-\alpha-1}^{-1}\frac{s-\polyr}{\alpha+1} \frac{\polyr-\alpha}{s-\alpha} \binom{\polyr-\alpha-1}{s-\alpha-1},\] this is equivalent to
\[\textstyle \binom{\alpha}{s-\alpha-1}^{-1}\frac{s-\polyr}{\alpha+1} \binom{\polyr-\alpha-1}{s-\alpha-1} + \sum_{j=1}^\alpha \frac{(-1)^{j+1}(\alpha-j+1)}{j+1} \binom{s-\alpha}{\alpha-j+1} \binom{\polyr-\alpha+j}{s-\alpha} = (-1)^{\alpha+1}. \]
The polynomial on the left side has degree at most \m{s-\alpha}. The coefficient of \m{\polyr^{s-\alpha}} in it is \m{-\frac{(2\alpha-s+1)!}{(\alpha+1)!}} plus \begin{align*} %
\textstyle %
\frac{1}{(s-\alpha-1)!} \sum_{j} \frac{(-1)^{j+1}}{j+1}\binom{s-\alpha-1}{\alpha-j} %
&  \textstyle\tempheightaer = %
\frac{1}{(s-\alpha-1)!} \sum_{j} \frac{(-1)^{j+1}}{j+1}[X^{s-2\alpha+j-1}] (1+X)^{s-\alpha-1}%
\\ &  \textstyle\tempheightaer = %
 \sum_{j} \frac{(-1)^{j+1}}{(j+1)}[X^{j}] X^{2\alpha-s-1}(1+X)^{s-\alpha-1}
\\ &  \textstyle\tempheightaer = %
 \frac{1}{(s-\alpha-1)!} \int_0^{-1} Y^{2\alpha-s+1}(1+Y)^{s-\alpha-1} \,\ddd Y \\ &  \textstyle\tempheightaer  = \frac{(-1)^{s} (2\alpha-s+1)!}{(\alpha+1)!}.
\end{align*} Since \m{s} is even,  that coefficient is zero. %
 Therefore it is enough to show that the two polynomials are equal when evaluated at \m{\polyr \in \{\alpha+1,\ldots,s\}}. %
  At these points the polynomial on the left side is equal to
 \[\textstyle %
 (s-\alpha) \sum_{j} \frac{(-1)^{j+1}}{j+1} \binom{s-\alpha-1}{\alpha-j} \binom{s-\alpha-\gamma+j}{s-\alpha}\]
 for \m{\gamma \in \{0,\ldots,s-\alpha-1\}}. We have
 \begin{align*}
     \textstyle \sum_{j} \frac{(-1)^{j+1}}{j+1} \binom{s-\alpha-1}{\alpha-j} \binom{s-\alpha-\gamma+j}{s-\alpha} & \textstyle = \sum_{j} \frac{(-1)^{s-\alpha+j+1}}{j+1} \binom{s-\alpha-1}{\alpha-j} \binom{\gamma-j-1}{s-\alpha}
     \\& \textstyle =  \sum_{u} \binom{\gamma}{u} \sum_{j} \frac{(-1)^{s-\alpha+j+1}}{j+1} \binom{s-\alpha-1}{\alpha-j} \binom{-j-1}{s-\alpha-u}
           \\& \textstyle = \sum_{u} \binom{\gamma}{u} \sum_{j} \frac{(-1)^{s-\alpha+j}}{s-\alpha-u} \binom{s-\alpha-1}{\alpha-j} \binom{-j-2}{s-\alpha-u-1}
                      \\& \textstyle = \sum_{u} \frac{(-1)^{u+1}}{s-\alpha-u}  \binom{\gamma}{u} \sum_{j} \binom{s-\alpha-1}{\alpha-j} \binom{-s+\alpha+u}{j+2}
                                            \\& \textstyle = \sum_{u} \frac{(-1)^{u+1}}{s-\alpha-u}  \binom{\gamma}{u} \binom{u-1}{\alpha+2}
                                                     =  \frac{(-1)^{\alpha+1}}{s-\alpha}.
 \end{align*}
The third equality follows from \m{\binom{\gamma}{u} = 0} for \m{u>s-\alpha-1}, and the last equality follows from \m{\binom{u-1}{\alpha+2}=0} for \m{u \in \{1,\ldots,s-\alpha-1\}}. In particular, \eqref{eqn423} is indeed  true.

 So both \m{F_1(\polyr)} and \m{F_2(\polyr)} have degree at most \m{\alpha+1}, and therefore they are equal if they are equal when evaluated at \[\textstyle\polyr \in \{s+\gamma(p-1) \,|\, \gamma \in \{0,\ldots,\alpha+1\}\}.\] It is easy to verify that \m{F_1(s)=F_2(s)}, %
  and when \[\textstyle \polyr \in \{s+\gamma(p-1) \,|\, \gamma \in \{1,\ldots,\alpha+1\}\}\] the fact that
 \[\textstyle \sum_{i=1}^{\gamma-1} \binom{s+\gamma(p-1)-\alpha+j}{i(p-1)+j}\binom{i(p-1)}{w} = F_{w,j}(s+\gamma(p-1))\]
  (due to~\eqref{ceqn6}) implies that 
  the equation \m{F_1(s+\gamma(p-1)) = F_2(s+\gamma(p-1))} is equivalent to
\begin{align*} & \textstyle \sum_{j=0}^\alpha \polyconst_j(s+\gamma(p-1)) \sum_{i=1}^{\gamma-1}  \binom{s+\gamma(p-1)-\alpha+j}{i(p-1)+j}\binom{i(p-1)}{w} %
= -\VERIFYIF{w=0} \binom{-\gamma(p-1)-1}{\alpha+1}. %
\end{align*} Note that \m{\binom{-\gamma(p-1)-1}{\alpha+1} = \binom{\gamma-1}{\alpha+1}=0} %
 and therefore the right side vanishes. Let us reiterate that all computations done in this proof are over \m{\FF_p}.
  Let us write \m{C_j^\gamma = \polyconst_j(s+\gamma(p-1))}. 
  The desired identity 
  \begin{align*} & \textstyle \sum_{j=0}^\alpha C_j^\gamma \sum_{i=1}^{\gamma-1}  \binom{s+\gamma(p-1)-\alpha+j}{i(p-1)+j}\binom{i(p-1)}{w} %
= 0 %
\end{align*}
  follows if
 \[\textstyle  \sum_{j=0}^\alpha C_j^\gamma  \binom{s+\gamma(p-1)-\alpha+j}{i(p-1)+j} = 0 \]
 for all \m{i \in \{1,\ldots,\gamma-1\}}. 
 If \m{j>0} and \m{C_j^\gamma \ne 0} then \[\textstyle j\GEQSLANT 2\alpha-s+1 \GEQSLANT \alpha+\gamma-s\] and consequently
 \begin{align*}\textstyle \binom{s+\gamma(p-1)-\alpha+j}{i(p-1)+j} & \textstyle = \left\{\begin{array}{rl}
     \binom{\gamma}{i} \binom{s-\alpha-\gamma+j}{s-\alpha-\gamma+i} & \text{if } s-\alpha-\gamma+i\GEQSLANT 0 \\[2pt]
    \binom{\gamma}{i-1} \binom{s-\alpha-\gamma+j}{p+s-\alpha-\gamma+i} & \text{if }s-\alpha-\gamma+i< 0
\end{array}\right. \\[5pt] & \textstyle = \binom{\gamma}{i} \binom{s-\alpha-\gamma+j}{s-\alpha-\gamma+i}. %
 \end{align*}
 On the other hand,
 \[\textstyle \binom{s+\gamma(p-1)-\alpha}{i(p-1)} = \binom{\gamma-1}{i-1} \binom{s-\alpha-\gamma}{s-\alpha-\gamma+i}.\]
 Since \[\textstyle \binom{\gamma-1}{i-1} = \frac{i}{\gamma} \binom{\gamma}{i} \in \FF_p^\times\] (as that \m{0 < i < \gamma \LEQSLANT \alpha+1}), what we want to show is that
 \[\textstyle C_0^\gamma \frac{i}{\gamma}\binom{s-\alpha-\gamma}{s-\alpha-\gamma+i} + \sum_{j=1}^{\alpha} C_j^\gamma \binom{s-\alpha-\gamma+j}{s-\alpha-\gamma+i} = 0\]
 for all \m{i \in \{1,\ldots,\gamma-1\}}.  That is equivalent to \[\textstyle F_3(s+\gamma(p-1)) = 0,\] where \m{F_3(z)\in\FF_p[z]} is defined as
  \[\textstyle F_3(\polyr) = \binom{\alpha}{s-\alpha-1}^{-1}\frac{s-\polyr-w}{\alpha+1} \binom{\polyr-\alpha-1}{s-\alpha-w-1} + \sum_{j=1}^\alpha \frac{(-1)^{j+1}(s-\alpha-w)}{j+1} \binom{s-\alpha-1}{\alpha-j} \binom{\polyr-\alpha+j}{s-\alpha-w}\]
  with \m{w = \gamma-i > 0}. The degree of \m{F_3(\polyr)} is at most \m{s-\alpha-w}, and in fact the coefficient of \m{\polyr^{s-\alpha-w}} in it is \m{-\frac{(s-\alpha-1)_w(2\alpha-s+1)!}{(\alpha+1)!}} plus \begin{align*} %
\textstyle %
\frac{1}{(s-\alpha-w-1)!} \sum_{j} \frac{(-1)^{j+1}}{j+1}\binom{s-\alpha-1}{\alpha-j} =  %
 \frac{(s-\alpha-1)_w(2\alpha-s+1)!}{(\alpha+1)!},
\end{align*}
 \IE the coefficient of \m{\polyr^{s-\alpha-w}} in it is zero. Therefore the degree of \m{F_3(\polyr)} is less than \m{s-\alpha-w}, so it is enough to show that  \m{F_3(\polyr)} is equal to zero when evaluated at \[\textstyle \polyr \in \{\alpha+1,\ldots,s-w\}.\] At these points \m{F_3(\polyr)} is equal to
  \[\textstyle %
 (s-\alpha-w) \sum_{j} \frac{(-1)^{j+1}}{j+1} \binom{s-\alpha-1}{\alpha-j} \binom{s-\alpha-\gamma+j}{s-\alpha-w}\]
 for \m{\gamma \in \{w,\ldots,s-\alpha-1\}}. We have
 \begin{align*}
     \textstyle \sum_{j} \frac{(-1)^{j+1}}{j+1} \binom{s-\alpha-1}{\alpha-j} \binom{s-\alpha-\gamma+j}{s-\alpha-w} & \textstyle = \sum_{j} \frac{(-1)^{s-\alpha+j-w+1}}{j+1} \binom{s-\alpha-1}{\alpha-j} \binom{\gamma-j-w-1}{s-\alpha-w}
     \\& \textstyle =  \sum_{u} \binom{\gamma-w}{u} \sum_{j} \frac{(-1)^{s-\alpha+j-w+1}}{j+1} \binom{s-\alpha-1}{\alpha-j} \binom{-j-1}{s-\alpha-u-w}
           \\& \textstyle = \sum_{u} \binom{\gamma-w}{u} \sum_{j} \frac{(-1)^{s-\alpha+j-w}}{s-\alpha-u-w} \binom{s-\alpha-1}{\alpha-j} \binom{-j-2}{s-\alpha-u-w-1}
                      \\& \textstyle = \sum_{u} \frac{(-1)^{u+1}}{s-\alpha-u-w}  \binom{\gamma-w}{u} \sum_{j} \binom{s-\alpha-1}{\alpha-j} \binom{-s+\alpha+u+w}{j+2}
                                            \\& \textstyle = \sum_{u} \frac{(-1)^{u+1}}{s-\alpha-u-w}  \binom{\gamma-w}{u} \sum_{j} \binom{u+w-1}{\alpha+2}
                                 =  0.
 \end{align*}
 The last equality follows from \m{\binom{\gamma-w}{u} = 0} for \[\textstyle u \not \in \{0,\ldots,s-\alpha-w-1\}.\] This proves that indeed \m{F_3(\polyr) = 0} and therefore that \m{F_1(z) = F_2(z)}.
\end{proof}

\begin{lemma}\label{lemmaComb555}
Suppose that \m{\alpha\in\ZZ_{\GEQSLANT0}}. For \m{w,j \in \{0,\ldots,\alpha\}} let \[\textstyle F_{w,j}(\polyr,\polys) \in \FF_p[\polyr,\polys]\] denote the polynomial
\[\textstyle \sum_v (-1)^{w-v} \binom{j+w-v-1}{w-v}\binom{\polyr-\alpha+j}{v} \left(\binom{\polys-\alpha+j-v}{j-v}-\binom{\polyr-\alpha+j-v}{j-v}\right).\]
Note that this depends on \m{\alpha}. Then
\begin{align*} \textstyle
     \sum_{j=1}^\alpha (-1)^{\alpha-j} \binom{\polys-\alpha+1}{\alpha-j} F_{w,j}(\polyr,\polys)  = 
  (-1)^\alpha(\VERIFYIF{w=\alpha}-\VERIFYIF{w=0}) \binom{\polys-\polyr}{\alpha}.
\end{align*}
\end{lemma}

\begin{proof}{Proof. }
   Both sides of the equation we want to prove have degree \m{\alpha} and the coefficient of \m{\polyr^\alpha} on each side is \m{\frac{1}{\alpha!}\left(\VERIFYIF{w=\alpha}-\VERIFYIF{w=0}\right)}. So the two sides are equal if they are equal when evaluated at the points \m{(\polyr,\polys)} such that
\[\textstyle (\polyr,\polys) \in \left\{ (u+\gamma(p-1)+\alpha,u+\alpha) \,|\, u \in \{0,\ldots,\alpha\} , \, \gamma \in \{0,\ldots,\alpha-1\}\right\}.\]
    The right side is zero when evaluated at these points, and
    \[\textstyle F_{w,j} (u+\gamma(p-1)+\alpha,u+\alpha) = \sum_{i=1}^\gamma \binom{u+\gamma(p-1)+j}{i(p-1)+j} \binom{i(p-1)}{w}\]
    by~\eqref{ceqn6}. Thus we want to show that
    \[\textstyle \sum_{j=1}^\alpha (-1)^{\alpha-j} \binom{u+1}{\alpha-j} \sum_{i=1}^\gamma \binom{u+\gamma(p-1)+j}{i(p-1)+j}\binom{i(p-1)}{w} = \BIGO{p}\] for \m{0\LEQSLANT u,w\LEQSLANT \alpha} and \m{0\LEQSLANT \gamma< \alpha}. Since \[\textstyle\binom{u+\gamma(p-1)+j}{i(p-1)+j}\binom{i(p-1)}{w} = \binom{\gamma}{i}\binom{u+j-\gamma}{j-i}\binom{-i}{w}+\BIGO{p},\] that is equivalent to
\[\textstyle \sum_{i,j>0} (-1)^{\alpha+w-i} \binom{u+1}{\alpha-j} \binom{\gamma}{i}\binom{\gamma-u-i-1}{j-i}\binom{i+w-1}{w} = \BIGO{p}.\] 
This follows from the facts that
\[\textstyle \sum_{j>0} \binom{u+1}{\alpha-j}\binom{\gamma-u-i-1}{j-i} = \binom{\gamma-i}{\alpha-i}\]
for \m{i>0} by Vandermonde's convolution formula, and \[\textstyle  \binom{\gamma}{i} \binom{\gamma-i}{\alpha-i} = \binom{\alpha}{i} \binom{\gamma}{\alpha}=0\] since \m{\gamma\in\{0,\ldots,\alpha-1\}}.
\end{proof}

  \begin{lemma}\label{comblemma1p} 
Suppose that \m{s,\alpha,\beta\in\ZZ} are such that
\[\textstyle s \in \{2,4,\ldots,p-3\} \text{ and } \alpha=\beta=\frac s2+1.\]
Let \m{M} denote the \m{(\alpha+1)\times(\alpha+1)} matrix with entries in \m{\FF_p} defined in lemma~\ref{comblemma1}. Suppose that \m{C_0,\ldots,C_\alpha\in \FF_p} are defined as
\begin{align*} \textstyle
 C_j =  (-1)^{\alpha+j+1}\alpha\binom{\alpha-2}{j-2}.\end{align*}
 Then \[\textstyle M(C_0,\ldots,C_\alpha)^T = (0,\ldots,0,1)^T.\]
\end{lemma}

\begin{proof}{Proof. }
  The equation associated with the \m{i}th row of \m{M} is straightforward if \m{i \not\in\{0,\alpha\}}. 
   Since \m{M_{0,j}} is equal to
    \begin{align*} \textstyle  & \tempheightaer\textstyle  \sum_{l,v=0}^\alpha (-1)^{l+v} \binom{j-v}{l-v} \binom{j-2}{v}^\partial \binom{\MORE{\alpha+j-v-2}}{\LESS{j-v}} (1+\VERIFYIF{\alpha=2 \,\& \, j=v})
 \\ & \null\qquad \textstyle 
  -  \binom{j-2}{j}^\partial -  \binom{j-2}{\alpha-2} \binom{0}{\alpha}^\partial
 \end{align*}
 and since
 \begin{align*}
     \textstyle \sum_j (-1)^j \binom{\alpha-2}{j-2} \frac{1}{j(j-1)}  & \textstyle = \frac1\alpha, \\
     \textstyle  \sum_j (-1)^j \binom{\alpha-2}{j-2} \binom{j-2}{s-\alpha} \binom{0}{\alpha}^\partial & \textstyle = (-1)^\alpha \binom{0}{\alpha}^\partial = -\frac{1}{\alpha},
      \\
     \textstyle \sum_{l}^\alpha (-1)^{l+\alpha} \binom{2-\alpha}{l-\alpha} \binom{0}{\alpha}^\partial \binom{p-1}{\alpha-2} %
     & \textstyle = -\frac{1%
     }{\alpha},
     \\
          \textstyle \sum_{l}^\alpha (-1)^{l+\alpha} \binom{2-\alpha}{l-\alpha} \binom{0}{\alpha}^\partial \binom{0}{\alpha-2}
     & \textstyle = -\frac{\VERIFYIF{\alpha=2}}{\alpha},
 \end{align*}
    the equation associated with the zeroth row is
    \[\textstyle \sum_j (-1)^j \binom{\alpha-2}{j-2} \sum_{l,v=0}^\alpha (-1)^{l+v} \binom{j-v}{l-v} \binom{j-2}{v}^\partial \binom{\alpha+j-v-2}{j-v} %
= -  %
 \frac 1\alpha, %
 \]
   and it follows from the fact that \[\textstyle \sum_{l=0}^\alpha (-1)^{l} \binom{l}{i} \binom{j-v}{l-v}= (-1)^{j} \binom{v}{j-i}\] for \m{0\LEQSLANT v\LEQSLANT j\LEQSLANT\alpha}. This shows the equation associated with the zeroth row. Since \m{M_{\alpha,j}} is equal to
           \begin{align*} \textstyle  & \tempheightaer\textstyle  \VERIFYIF{j=2}   \binom{0}{\alpha}^\partial +  \VERIFYIF{j=\alpha} F_{\alpha,\alpha,0}(\alpha-2)
    - \binom{j-2}{\alpha-2}^\partial 
   - (-1)^{\alpha}  \binom{j-2}{\alpha-2} \binom{-1}{\alpha}^\partial, \end{align*}
the equation associated with the \m{\alpha}th row is
   \[\textstyle \binom{0}{\alpha}^\partial
    - \sum_j (-1)^j \binom{\alpha-2}{j-2}\binom{j-2}{\alpha-2}^\partial 
   = \binom{-1}{\alpha}^\partial - (-1)^\alpha  F_{\alpha,\alpha,0}(\alpha-2) + \frac{(-1)^{\alpha+1}}{\alpha},\]
   and it follows from the facts that
    \[\textstyle F_{\alpha,\alpha,0}(\alpha-2) =  F_{\alpha,\alpha,0}(-1) = (-1)^\alpha \binom{-1}{\alpha}^\partial\]
    and that the polynomial \m{\binom{X}{\alpha-2}^\partial\in \FF_p[X]} has degree less than \m{\alpha-2} (and is zero if \m{\alpha=2}) and therefore \[\textstyle \sum_j (-1)^j \binom{\alpha-2}{j-2}\binom{j-2}{\alpha-2}^\partial = 0.\]This shows the equation associated with the \m{\alpha}th row and concludes the proof.
\end{proof}

\begin{lemma}\label{comblemma1cc} 
Suppose that \m{s,\alpha,\beta\in\ZZ} are such that
\[\textstyle s \in \{2,4,\ldots,p-3\} \text{ and } \alpha=\frac s2+1 \text{ and } \beta \in \{\frac s2,\frac s2+1\}.\]
Let \m{A_0} denote the \m{\beta\times\beta} matrix with entries in \m{\QQ_p} defined as
\[\textstyle A_0 =  \left(p^{\VERIFYIF{j=1}-\VERIFYIF{i\LEQSLANT\beta-\alpha+1}}\binom{s+\beta(p-1)-\alpha+j}{i(p-1)+j}\right)_{0\LEQSLANT i< \beta, \, \alpha-\beta<j\LEQSLANT\alpha}.\]
Then \m{A_0} has entries in \m{\ZZ_p} and is invertible over \m{\ZZ_p}.
\end{lemma}

\begin{proof}{Proof. }
 It is easy to verify that \m{A_0} is integral, since if \m{j>1} then \[\textstyle s-\alpha-\beta+j \GEQSLANT 0\] and therefore
 \[\textstyle \binom{s+\beta(p-1)-\alpha+j}{i(p-1)+j} = \binom{\beta}{i} \binom{s-\alpha-\beta+j}{j-i} + \BIGO{p} = \BIGO{p}\]
 for \m{i \LEQSLANT \beta-\alpha+1}. Let us show that \m{A_0} is invertible (over \m{\ZZ_p}) by showing that \m{\overline A_0} is invertible (over \m{\FF_p}). Suppose first that \m{\beta=\alpha-1} and denote the columns of \m{A_0} by \m{\colA_2,\ldots,\colA_\alpha}. The bottom left \m{(\alpha-3)\times(\alpha-3)} submatrix of \m{\overline A_0} is upper triangular with units on the diagonal. Moreover, since
  \begin{align*}
 \textstyle\sum_j (-1)^j \binom{s-\beta-j-1}{\alpha-i-j-1} \binom{\alpha-2}{j}     &\textstyle= \sum_j (-1)^j \binom{\beta-j-1}{i-1} \binom{\beta-1}{j} = 0,\\
     \textstyle \sum_j (-1)^{j-1} (j-1) \binom{s-\beta-j}{\alpha-i-j} \binom{\alpha-1}{j}    & 
     \textstyle =      \textstyle \sum_j (-1)^{j-1} (j-1) \binom{\beta-j}{i-1} \binom{\beta}{j} = 0,
 \end{align*}
  all but the top two entries of each of the vectors
 \begin{align*}
    & \textstyle \colA_{\alpha-1} - \binom{\alpha-2}{1} \colA_{\alpha-2} + \cdots + (-1)^{\alpha-3}\binom{\alpha-2}{\alpha-3}\colA_2,\\
     & 
     \textstyle \colA_{\alpha} - \binom{\alpha-1}{2} \colA_{\alpha-2} + \cdots + (-1)^{\alpha-3} (\alpha-3) \binom{\alpha-1}{\alpha-2}\colA_2
 \end{align*}
 are zero. Thus it is enough to show that the \m{2\times 2} matrix consisting of those four entries is invertible (over \m{\FF_p}). This \m{2\times 2} matrix is
\[\textstyle  \left(\begin{smallmatrix} e_{0,0}  & e_{0,1} \\  (-1)^{\beta}\beta &  (-1)^{\beta}\beta(\beta-1)  \end{smallmatrix}\right)\]
with
\begin{align*}
    \textstyle e_{0,0} & \textstyle = \beta\sum_{j=0}^{\beta-1} (-1)^j \binom{\beta-1}{j} \binom{\beta-j-1}{\beta-j}^\partial = \sum_{j=0}^{\beta-1} \frac{(-1)^j\beta}{\beta-j} \binom{\beta-1}{j} \\ & \textstyle\null\qquad = \sum_{j=0}^{\beta-1} (-1)^j \binom{\beta}{j} = (-1)^{\beta+1},\\[5pt]
    \textstyle e_{0,1} & \textstyle = \beta \sum_{j=0}^{\beta} (-1)^{j-1} (j-1) \binom{\beta}{j} \binom{\beta-j}{\beta-j+1}^\partial = \sum_{j=0}^{\beta} \frac{(-1)^{j-1} \beta (j-1)}{\beta-j+1} \binom{\beta}{j}  \\ & \textstyle\null\qquad = \sum_{j=0}^{\beta} \frac{(-1)^{j-1} \beta (j-1)}{\beta+1} \binom{\beta+1}{j} = \frac{(-1)^{\beta-1}\beta^2}{\beta+1},
\end{align*}
so it has determinant \m{\frac{\beta}{\beta+1} \in \FF_p^\times}. Now suppose that \m{\beta=\alpha} and denote the columns of \m{A_0} by \m{\colA_1,\ldots,\colA_\alpha}. The bottom left \m{(\alpha-1)\times(\alpha-1)} submatrix of \m{\overline A_0} is upper triangular with units on the diagonal,
  all but the top entry of the vector
 \begin{align*}
    & \textstyle \colA_{\alpha} - \binom{\alpha-2}{1} \colA_{\alpha-1} + \cdots + (-1)^{\alpha-2} \binom{\alpha-2}{\alpha-2}\colA_2
 \end{align*}
 are zero, and that top entry is
 \[\textstyle  \beta\sum_{j=0}^{\beta-2} (-1)^j \binom{\beta-2}{j} \binom{\beta-j-2}{\beta-j}^\partial = \sum_{j=0}^{\beta-2} \frac{(-1)^j}{\beta-1} \binom{\beta}{j} %
 = (-1)^\beta \in \FF_p^\times.\] Therefore \m{\overline A_0} is  invertible.
\end{proof}

\begin{lemma}\label{comblemma1} 
Suppose that \m{s,\alpha,\beta\in\ZZ} are such that
\[\textstyle s \in \{2,4,\ldots,p-3\} \text{ and } \frac s2 \LEQSLANT \alpha \LEQSLANT s \text{ and } 1\LEQSLANT \beta \LEQSLANT \alpha.\]
Let \m{M} denote the \m{(\alpha+1)\times(\alpha+1)} matrix with entries in \m{\FF_p} such that if \m{i\in\{1,\ldots,\beta-1\}} and \m{j\in\{0,\ldots,\alpha\}} then
\[\textstyle M_{i,j} = \binom{s+\beta(p-1)-\alpha+j}{i(p-1)+j},\]
and if \m{i\in\{0,\ldots,\alpha\}\backslash\{1,\ldots,\beta-1\}} and \m{j\in\{0,\ldots,\alpha\}} then
\begin{align*} \textstyle M_{i,j} & \tempheightaer\textstyle = \sum_{l,v=0}^\alpha (-1)^{i+l+v} \binom{l}{i} \binom{j-v}{l-v} \binom{s-\alpha-\beta+j}{v}^\partial \binom{\MORE{s-\alpha+j-v}}{\LESS{j-v}} %
 \\ & \tempheightaer \textstyle\null\qquad - \VERIFYIF{i=0} \binom{s-\alpha-\beta+j}{j}^\partial - \VERIFYIF{i=\beta} \binom{s-\alpha-\beta+j}{s-\alpha}^\partial \\ & \tempheightaer\textstyle\null\qquad
  - (-1)^{i}  \binom{s-\alpha-\beta+j}{s-\alpha} \sum_{l=0}^\alpha \binom{l}{i} \binom{l-\beta-1}{l}^\partial.\end{align*}
Suppose that \m{C_0,\ldots,C_\alpha \in \FF_p} are defined as
\begin{align*} 
 C_j = \left\{\begin{array}{rl} 1
& \text{if } j=0,
\\[5pt] \frac{(-1)^{j+1}(s-\alpha-\beta)}{\beta} \binom{j}{2\alpha-s+1} \binom{\alpha+1}{j+1} & \text{if } j \in \{1,\ldots,\alpha\}. \end{array}\right.\end{align*}
Then
\[\textstyle M (C_{0},C_1,\ldots,C_\alpha)^T = 
\left(\frac{(-1)^{\alpha+\beta+1}(s-\alpha)(\alpha-\beta+1)}{\beta^2(2\alpha-s+1) \binom{\alpha}{\beta}} \binom{\alpha}{s-\alpha},0,\ldots,0\right)^T.
\] 
\end{lemma}

\begin{proof}{Proof. }
Let us denote the rows of \m{M} by \[\textstyle \rowR_0,\ldots,\rowR_\alpha.\] Note that if \m{j>0} and \m{C_j \ne 0} then \m{j>2\alpha-s}, so \m{\textstyle s-\alpha+j > \alpha} and in particular
  \[\textstyle \binom{\MORE{s-\alpha+j-v}}{\LESS{j-v}} = \binom{s-\alpha+j-v}{j-v}.\]
   We have the following string of equations:
 \begin{align*}
     & \tempheightaer \textstyle\sum_{j\GEQSLANT 0} (-1)^{j+1} \binom{j}{2\alpha-s+1} \binom{\alpha+1}{j+1}  \binom{s-\alpha-\beta+j}{s-\alpha}
     \\ & \tempheightaer\textstyle\null\qquad = \sum_u \binom{\beta }{u} \sum_{j\GEQSLANT 0} (-1)^{s-\alpha+j+1} \binom{j}{2\alpha-s+1} \binom{\alpha+1}{j+1}  \binom{-j-1}{s-\alpha-u}
          \\ & \tempheightaer\textstyle\null\qquad = \sum_u \binom{\beta }{u} \binom{\alpha-u+1}{s-\alpha-u} \sum_{j\GEQSLANT 0} (-1)^{j+u+1} \binom{\alpha+1}{j+1}  \binom{s-\alpha+j-u}{\alpha-u+1} 
           \\ & \tempheightaer\textstyle\null\qquad = \sum_u \binom{\beta }{u} \binom{\alpha-u+1}{s-\alpha-u} \left((-1)^{u+1}\binom{s-\alpha-u-1}{\alpha-u+1} + \sum_{j} (-1)^{j+u+1} \binom{\alpha+1}{j+1}  \binom{s-\alpha+j-u}{\alpha-u+1} \right)
                      \\ & \tempheightaer\textstyle\null\qquad = \sum_u \binom{\beta }{u} \binom{\alpha-u+1}{s-\alpha-u} \left((-1)^{u+1}\binom{s-\alpha-u-1}{\alpha-u+1} + \VERIFYIF{u=0}(-1)^{\alpha+1} \right)
                       \\ & \tempheightaer\textstyle\null\qquad = (-1)^\alpha\left(\binom{\beta}{s-\alpha}-\binom{\alpha+1}{s-\alpha}\right).
  \end{align*}
  The first two equalities amount to rewriting the binomial coefficients. The third equality amounts to computing the inner sum. The fourth equality follows from~\eqref{ceqn5}. The fifth equality amounts to computing the outer sum.
 This string of equations implies that
  \[\textstyle \sum_{j=0}^\alpha C_j \binom{s-\alpha-\beta+j}{s-\alpha} = (-1)^{\alpha+1} \frac{s-\alpha-\beta}{\beta} \binom{\alpha+1}{s-\alpha}.\]
  Our task is to compute \m{\rowR_i (C_{0},C_1,\ldots,C_\alpha)^T} for \m{i \in\{0,\ldots,\alpha\}}.
  
  \begin{itemize}[leftmargin=*]
      \item \emph{Computing \m{\rowR_0 (C_{0},C_1,\ldots,C_\alpha)^T}.}   If \m{j>2\alpha-s} then \[\textstyle \sum_{l=0}^\alpha (-1)^l \binom{l}{i} \binom{j-v}{l-v} = (-1)^j \binom{v}{j-i}\] for \m{0 \LEQSLANT v \LEQSLANT j \LEQSLANT \alpha} and therefore
\begin{align*} \textstyle M_{0,j} & \tempheightaer\textstyle = \sum_{v=0}^\alpha (-1)^{j+v} \binom{v}{j} \binom{s-\alpha-\beta+j}{v}^\partial \binom{{s-\alpha+j-v}}{{j-v}} %
 \\ & \tempheightaer \textstyle\null\qquad -  \binom{s-\alpha-\beta+j}{j}^\partial -  \binom{s-\alpha-\beta+j}{s-\alpha} \sum_{l=0}^\alpha  \binom{l-\beta-1}{l}^\partial
 \\ & \tempheightaer\textstyle =   - \binom{s-\alpha-\beta+j}{s-\alpha} \sum_{l=0}^\alpha \binom{l-\beta-1}{l}^\partial.
 \end{align*}
 The second equality follows from the fact that \m{\binom{v}{j} = 0} if \m{v < j}. We also have
 \begin{align*} \textstyle M_{0,0} & \tempheightaer\textstyle =  \sum_{l,v=0}^\alpha (-1)^{l+v} \binom{-v}{l-v} \binom{s-\alpha-\beta}{v}^\partial \binom{\MORE{s-\alpha-v}}{\LESS{-v}}
 \\ & \tempheightaer \textstyle\null\qquad -  \binom{s-\alpha-\beta}{s-\alpha} \sum_{l=0}^\alpha \binom{l-\beta-1}{l}^\partial
 \\ & \tempheightaer\textstyle =   \sum_{l,v=0}^\alpha (-1)^{\alpha+l+v} \binom{-v}{l-v} \binom{s-\alpha-\beta}{v}^\partial \binom{v}{s-\alpha}
 \\ & \tempheightaer \textstyle\null\qquad - \binom{s-\alpha-\beta}{s-\alpha} \sum_{l=0}^\alpha \binom{l-\beta-1}{l}^\partial
  \\ & \tempheightaer\textstyle =   \sum_{l,v=0}^\alpha (-1)^{\alpha+l+v} \binom{-v}{l-v} \binom{s-\alpha-\beta}{v}^\partial \binom{v}{s-\alpha}
 \\ & \tempheightaer \textstyle\null\qquad +  \binom{s-\alpha-\beta}{s-\alpha} \sum_{l=0}^\alpha (-1)^l \sum_v (-1)^{v+1}\binom{s-\alpha-\beta}{v}^\partial \binom{s-\alpha-v}{l-v}.
 \end{align*} The third equality follows from lemma~\ref{part2X3a}.
  Thus \m{\rowR_0 (C_0,\ldots,C_\alpha)^T} is equal to
 \begin{align*} & \textstyle  \sum_{l,v=0}^\alpha (-1)^{\alpha+l+v} \binom{s-\alpha-\beta}{v}^\partial \left(\binom{-v}{l-v}  \binom{v}{s-\alpha} + \frac{s-\alpha-\beta}{\beta} \binom{\alpha+1}{s-\alpha} \binom{s-\alpha-v}{l-v}\right)
  \\ & \textstyle\null\qquad\tempheightaer = (-1)^\alpha \sum_{v=0}^\alpha \binom{s-\alpha-\beta}{v}^\partial \left(\binom{\alpha}{v} \binom{v}{s-\alpha} + \frac{s-\alpha-\beta}{\beta} \binom{\alpha+1}{s-\alpha} \binom{2\alpha-s}{\alpha-v}\right)
    \\ & \textstyle\null\qquad\tempheightaer = \left(     \binom{\alpha}{s-\alpha} + \frac{s-\alpha-\beta}{\beta} \binom{\alpha+1}{s-\alpha}  \right) \sum_{v=0}^\alpha (-1)^v \binom{s-\alpha-\beta}{v}^\partial \binom{s-\alpha-v-1}{\alpha-v}
                \\ & \textstyle\null\qquad\tempheightaer = \frac{(-1)^{\alpha+\beta+1}}{\beta \binom{\alpha}{\beta}}\left(     \binom{\alpha}{s-\alpha} + \frac{s-\alpha-\beta}{\beta} \binom{\alpha+1}{s-\alpha}  \right) 
                \\ & \textstyle\null\qquad\tempheightaer =   \frac{(-1)^{\alpha+\beta+1}(s-\alpha)(\alpha-\beta+1)}{\beta^2(2\alpha-s+1) \binom{\alpha}{\beta}} \binom{\alpha}{s-\alpha}.
\end{align*}
The third equality follows from lemma~\ref{part2X3a}. Thus we have computed
\[\textstyle \rowR_0 (C_0,\ldots,C_\alpha)^T = \frac{(-1)^{\alpha+\beta+1}(s-\alpha)(\alpha-\beta+1)}{\beta^2(2\alpha-s+1) \binom{\alpha}{\beta}} \binom{\alpha}{s-\alpha}.\]
\item \emph{Computing \m{\rowR_i (C_{0},C_1,\ldots,C_\alpha)^T} for \m{i \in \{1,\ldots,\beta-1\}}.} Let \m{w \in \ZZ} be such that \m{i = \beta-w \in \{1,\ldots,\beta-1\}}. Then
\begin{align*}
    & \textstyle \sum_{j\GEQSLANT 0} \frac{(-1)^{j+1}(s-\alpha-\beta)}{\beta} \binom{j}{2\alpha-s+1} \binom{\alpha+1}{j+1} \binom{s-\alpha-\beta+j}{s-\alpha-w}  
    \\ & \textstyle\null\qquad = \sum_u \binom{\alpha-\beta+1}{u} \sum_{j\GEQSLANT 0} \frac{(-1)^{j+1}(s-\alpha-\beta)}{\beta} \binom{j}{2\alpha-s+1} \binom{\alpha+1}{j+1} \binom{s-2\alpha+j-1}{s-\alpha-w-u} 
        \\ & \textstyle\null\qquad = \sum_u \binom{\alpha-\beta+1}{u} \binom{\alpha-w-u+1}{s-\alpha-w-u} \sum_{j\GEQSLANT 0} \frac{(-1)^{j+1}(s-\alpha-\beta)}{\beta} \binom{j}{\alpha-w-u+1} \binom{\alpha+1}{j+1}
                \\ & \textstyle\null\qquad = (-1)^{\alpha-w-u} \frac{s-\alpha-\beta}{\beta} \sum_u  \binom{\alpha-\beta+1}{u} \binom{\alpha-w-u+1}{s-\alpha-w-u}   
                             \\ & \textstyle\null\qquad = \frac{s-\alpha-\beta}{\beta} \sum_u  \binom{\alpha-\beta+1}{u} \binom{s-2\alpha-2}{s-\alpha-w-u}    
                          \\ & \textstyle\null\qquad = \frac{s-\alpha-\beta}{\beta}   \binom{s-\alpha-\beta-1}{s-\alpha-w} = -\frac{i}{\beta} \binom{s-\alpha-\beta}{s-\alpha-w}. 
\end{align*}
The third equality follows from~\eqref{ceqn5}. Consequently, if \m{i \in \{1,\ldots,\beta-1\}} then
\[\textstyle \rowR_i (C_0,\ldots,C_\alpha)^T = 0.\]
\item \emph{Computing \m{\rowR_i (C_{0},C_1,\ldots,C_\alpha)^T} for \m{i \in \{\beta,\ldots,\alpha\}}.} For these \m{i} we have
\begin{align*} \textstyle M_{i,0} & \tempheightaer\textstyle = \sum_{l,v=0}^\alpha (-1)^{\alpha+i+l+v} \binom{l}{i} \binom{v}{s-\alpha} \binom{s-\alpha-\beta}{v}^\partial \binom{-v}{l-v} %
 \\ & \tempheightaer \textstyle\null\qquad  - \VERIFYIF{i=\beta} \binom{s-\alpha-\beta}{s-\alpha}^\partial 
  - (-1)^{i} \binom{s-\alpha-\beta}{s-\alpha} \sum_{l=0}^\alpha  \binom{l}{i} \binom{l-\beta-1}{l}^\partial,
  \end{align*} 
  and for \m{j > 2\alpha-s} we also have
  \begin{align*} \textstyle M_{i,j} & \tempheightaer\textstyle = \sum_{l,v=0}^\alpha (-1)^{i+l+v} \binom{l}{i} \binom{j-v}{l-v} \binom{s-\alpha-\beta+j}{v}^\partial \binom{{s-\alpha+j-v}}{{j-v}} %
\\ & \tempheightaer\textstyle\null\qquad
   - \VERIFYIF{i=\beta} \binom{s-\alpha-\beta+j}{s-\alpha}^\partial - (-1)^{i} \binom{s-\alpha-\beta+j}{s-\alpha} \sum_{l=0}^\alpha \binom{l}{i} \binom{l-\beta-1}{l}^\partial.  \end{align*}
       The identity
  \begin{align*}
  &  \textstyle \sum_{j=0}^\alpha (-1)^{j+1} \binom{j}{2\alpha-s+1} \binom{\alpha+1}{j+1} \binom{\polyr+s-\alpha+j}{s-\alpha}^\partial
  \\ & \null\qquad \textstyle\tempheightaer = \frac{\partial}{\partial \polyr} \left(\sum_{j=0}^\alpha (-1)^{j+1} \binom{j}{2\alpha-s+1} \binom{\alpha+1}{j+1} \binom{\polyr+s-\alpha+j}{s-\alpha}\right)
    \\ & \null\qquad \textstyle\tempheightaer = \frac{\partial}{\partial \polyr} \left(\binom{\polyr+s-\alpha-1}{s-\alpha}-\binom{s-2\alpha-2}{s-\alpha}\right)
     \\ & \null\qquad \textstyle\tempheightaer =  \binom{\polyr+s-\alpha-1}{s-\alpha}^\partial
 \end{align*}
 is true over \m{\QQ_p[\polyr]}.
  By evaluating at \m{\polyr=-\beta} we get
  \begin{align*}%
  \textstyle  \sum_{j=1}^\alpha C_j \binom{s-\alpha-\beta+j}{s-\alpha}^\partial = \frac{s-\alpha-\beta}{\beta} \binom{s-\alpha-\beta-1}{s-\alpha}^\partial,\end{align*}
  and consequently
  \begin{align*}\textstyle \binom{s-\alpha-\beta}{s-\alpha}^\partial + \sum_{j=1}^\alpha C_j \binom{s-\alpha-\beta+j}{s-\alpha}^\partial = \frac{1}{\beta}\binom{s-\alpha-\beta-1}{s-\alpha-1}. \end{align*}
  This means that \m{(-1)^{\alpha+i}\beta \rowR_i(C_0,\ldots,C_\alpha)^T} is equal to \m{\polyH(-\beta)}, with
\begin{align*} \textstyle  & \tempheightaer\textstyle 
 \polyH(\polyr) = (\alpha-s) \polyH_1'(\polyr) - \polyH_2(\polyr) + (\polyr+s-\alpha) (\polyH_1'(\polyr) + \polyH_3'(\polyr)+  \polyH_4'(\polyr)) %
 \end{align*}
 and
  \begin{align*}
  \polyH_1(\polyr) &\tempheightaer \textstyle = \sum_{l,v=0}^\alpha (-1)^{l+v+1} \binom{l}{i} \binom{v}{s-\alpha} \binom{\polyr + s-\alpha}{v} \binom{-v}{l-v}, \\
  \polyH_2(\polyr)  &\tempheightaer \textstyle =  \binom{i-1}{s-\alpha-1}\binom{\polyr +i-1}{i-1}, \\
      \polyH_3(\polyr) &\tempheightaer \textstyle = \sum_{l,j,v=0}^\alpha   (-1)^{\alpha+j+l+v+1} \binom{j}{2\alpha-s+1}\binom{\alpha+1}{j+1} \binom{l}{i} \binom{j-v}{l-v} \binom{\polyr+s-\alpha+j}{v} \binom{{s-\alpha+j-v}}{{j-v}} , \\
      \polyH_4(\polyr) &\tempheightaer =\textstyle \binom{\alpha+1}{s-\alpha} \sum_{l=0}^\alpha \binom{l}{i} \binom{\polyr+l-1}{l} = \binom{\alpha+1}{s-\alpha} \binom{\polyr + \alpha}{\alpha-i} \binom{z+i-1}{i}.
  \end{align*}
  So we want to show that \m{\polyH(-\beta) = 0}. If \m{s=\alpha+\beta} then this equation amounts to
  \[\textstyle \beta \polyH_1'(-\beta)+\polyH_2(-\beta) = 0,\] and indeed
  \begin{align*}
    \beta \polyH_1'(-\beta)   & \textstyle = \beta \sum_{l,v=0}^\alpha (-1)^{l+v+1} \binom{l}{i} \binom{v}{\beta} \binom{0}{v}^\partial \binom{-v}{l-v}
    \\
     & \textstyle = \sum_{l,v=0}^\alpha \frac{(-1)^{l}\beta}{v} \binom{l}{i} \binom{v}{\beta}  \binom{-v}{l-v}\phantom{^\partial}
         \\
     & \textstyle = \sum_{l,v=0}^\alpha (-1)^{l} \binom{l}{i} \binom{v-1}{\beta-1}  \binom{-v}{l-v}\phantom{^\partial}
              \\
     & \textstyle = \sum_{l,v=0}^\alpha (\VERIFYIF{l=0} (-1)^{\beta+l+1} + \VERIFYIF{l=\beta}(-1)^l) \binom{l}{i}\phantom{\binom{\beta}{i}^\partial}   \\
     & \textstyle =  (-1)^\beta \binom{\beta}{i}\phantom{\binom{\beta}{i}^\partial}  \\
     & \textstyle = \VERIFYIF{i=\beta} (-1)^\beta\phantom{\binom{\beta}{i}^\partial}  \\
     & \textstyle = - \binom{i-1}{\beta-1} \binom{i-\beta-1}{i-1} = - \polyH_2(-\beta) .\phantom{\binom{\beta}{i}^\partial} \phantom{^\partial}
  \end{align*}
    Now suppose that \m{s\ne \alpha+\beta}. As in the proof of lemma~\ref{part2X3a} we can simplify \m{\polyH_1(\polyr)} to
  \begin{align*}
      \polyH_1(\polyr) & \textstyle = -\binom{\polyr+s-\alpha}{s-\alpha} \sum_{l=0}^\alpha \binom{l}{i} \binom{\polyr+l-1}{l+\alpha-s}.
  \end{align*}
    We can also simplify \m{\polyH_3(\polyr)} to
  \begin{align*}
      \polyH_3(\polyr) & \textstyle = \sum_{j,v=0}^\alpha   (-1)^{\alpha+v+1} \binom{j}{2\alpha-s+1}\binom{\alpha+1}{j+1} \binom{v}{j-i} \binom{\polyr+s-\alpha+j}{v} \binom{{s-\alpha+j-v}}{{j-v}}
      \\
      & \textstyle = \binom{{\polyr+i-1}}{{i}} \sum_{j=0}^\alpha   (-1)^{\alpha+j+1} \binom{j}{2\alpha-s+1}\binom{\alpha+1}{j+1} \binom{\polyr+s-\alpha+j}{j-i}.
  \end{align*}
  Suppose first that \m{i>\beta}. Then
    \begin{align*}
  \polyH_1'(-\beta) &\tempheightaer \textstyle = \smash{- \sum_{l=0}^\alpha \binom{l}{i}\Big(\binom{s-\alpha-\beta}{s-\alpha}^\partial \binom{l-\beta-1}{l+\alpha-s}+\binom{s-\alpha-\beta}{s-\alpha} \binom{l-\beta-1}{l+\alpha-s}^\partial\Big), } \\
    \polyH_2(-\beta) &\tempheightaer \textstyle = 0, \\
      \polyH_3'(-\beta) &\tempheightaer \textstyle =
    \binom{{i-\beta-1}}{{i}}^\partial \sum_{j=0}^\alpha   (-1)^{\alpha+j+1} \binom{j}{2\alpha-s+1}\binom{\alpha+1}{j+1} \binom{s-\alpha-\beta+j}{j-i}, \\
      \polyH_4'(-\beta) &\tempheightaer =\textstyle  \binom{\alpha+1}{s-\alpha} \binom{ \alpha-\beta}{\alpha-i} \binom{i-\beta-1}{i}^\partial.
  \end{align*}
  Thus if \m{s>\alpha+\beta} then the equation \m{\polyH(-\beta)=0} is equivalent to
 \[\textstyle L_1(s,\alpha,\beta,i) = R_1 (s,\alpha,\beta,i)\]
 with
 \begin{align*}
     \textstyle L_1 & \textstyle \defeq  \sum_{l=0}^\alpha \binom{l}{\beta+1} \binom{l-\beta-1}{i-\beta-1}  \binom{l-\beta-1}{s-\alpha-\beta-1}, \\
    R_1  & \textstyle \defeq \binom{s-\alpha}{\beta+1} \sum_{j=0}^\alpha   (-1)^{\alpha+j+1}  \binom{j}{2\alpha-s+1}\binom{\alpha+1}{j+1} \binom{s-\alpha-\beta+j}{j-i} \\ & \textstyle \null\qquad + \binom{\alpha+1}{\beta+1} \binom{\alpha-\beta}{s-\alpha-\beta-1} \binom{ \alpha-\beta}{\alpha-i}.
 \end{align*}
 Let us in fact show that \[\textstyle L_1(\zs,\za,\zb,\zi) = R_1(\zs,\za,\zb,\zi)\] for all \m{\zs,\za,\zb,\zi\GEQSLANT 0}. We clearly have \[\textstyle L_1(\zs,0,\zb,\zi) = R_1(\zs,0,\zb,\zi)\] since both sides are zero, and
  \begin{align*} & \textstyle R_1(\zs+1,\za+1,\zb,\zi) - R_1(\zs,\za,\zb,\zi) \\ & \textstyle \null\qquad\qquad - L_1(\zs+1,\za+1,\zb,\zi) + L_1(\zs,\za,\zb,\zi)
    \\ & \textstyle \null\qquad = \binom{\zs-\za}{\zb+1} \frac{\zs-\za}{2\za-\zs+2} \sum_{j=0}^{\za} (-1)^{\za+j} \binom{j}{2\za-\zs+1}\binom{\za+1}{j}\binom{\zs-\za-\zb+j}{j-\zi}
    \\ & \textstyle \null\qquad\qquad - \binom{\zs-\za}{\zb+1} \binom{\za+1}{2\za-\zs+2}\binom{\zs-\zb+1}{\za-\zi+1}
  \\ & \textstyle \null\qquad\qquad + \frac{\zs-\za}{2\za-\zs+2} \binom{\za+1}{\zb+1}\binom{\za-\zb}{\zs-\za-\zb-1}\binom{\za-\zb+1}{\zi-\zb}.
  \end{align*}
All we need to show is that this is zero for all \m{\zs,\za,\zb,\zi\GEQSLANT 0}, which follows from
\begin{align}\label{eq29835798}
\nonumber &  \textstyle \sum_{j} (-1)^j  \binom{j}{2\za-\zs+1}\binom{\za+1}{j} \binom{\zs-\za-\zb+j}{\zs-\za-\zb+\zi} \\ \nonumber &
  \textstyle 
  \null\qquad =  \textstyle \sum_{j,\varu} (-1)^j \binom{\za-\zb+1}{\zs-\za-\zb+i-\varu} \binom{j}{2\za-\zs+1}\binom{\za+1}{j} \binom{\zs-2\za+j-1}{\varu} 
   \\\nonumber & \textstyle
  \null\qquad =  \textstyle \sum_{j,\varu} (-1)^j \binom{\za-\zb+1}{\zs-\za-\zb+i-\varu} \binom{j}{2\za-\zs+\varu+1}\binom{\za+1}{j} \binom{2\za-\zs+\varu+1}{\varu}
     \\ \nonumber& \textstyle
  \null\qquad =  \textstyle \sum_{j,\varu} (-1)^{\zs+\varu+1} \binom{\za-\zb+1}{\zs-\za-\zb+\zi-\varu} \binom{\zs-2\za-\varu-2}{j+\zs-2\za-\varu-1}\binom{\za+1}{\za-j+1} \binom{2\za-\zs+\varu+1}{\varu}
       \\ \nonumber& \textstyle
  \null\qquad =  \textstyle \sum_{\varu} (-1)^{\zs+\varu+1} \binom{\za-\zb+1}{\zs-\za-\zb+\zi-\varu} \binom{\zs-\za-\varu-1}{\zs-\za-\varu} \binom{2\za-\zs+\varu+1}{\varu}
         \\ & \textstyle
  \null\qquad =  \textstyle (-1)^{\za+1} \binom{\za-\zb+1}{\zi-\zb} \binom{\za+1}{\zs-\za}.%
\end{align}
Similarly, if \m{s<\alpha+\beta} then the equation \m{\polyH(-\beta)=0} is equivalent to
 \[\textstyle L_2(s,\alpha,\beta,i) = R_2 (s,\alpha,\beta,i)\]
 with
 \begin{align*}
     \textstyle L_2 & \textstyle \defeq \sum_{l=\beta+1}^\alpha \binom{l}{s-\alpha} \binom{l-\beta-1}{i-\beta-1}, \\
    R_2  & \textstyle \defeq\sum_{j=0}^\alpha   (-1)^{\alpha+j+1}  \binom{j}{2\alpha-s+1}\binom{\alpha+1}{j+1} \binom{s-\alpha-\beta+j}{j-i}+  \binom{\alpha+1}{s-\alpha}\binom{ \alpha-\beta}{\alpha-i}.
 \end{align*}
  Let us in fact show that \[\textstyle L_2(\zs,\za,\zb,\zi) = R_2(\zs,\za,\zb,\zi)\] for all \m{\zs\GEQSLANT\za\GEQSLANT\zi>\zb\GEQSLANT 0}.
  It is easy to verify that \[\textstyle L_2(\zs,\zi,\zb,\zi) = R_2(\zs,\zi,\zb,\zi),\] and
    \begin{align*} & \textstyle R_2(\zs+1,\za+1,\zb,\zi) - R_2(\zs,\za,\zb,\zi) \\ & \textstyle\null\qquad\qquad - L_2(\zs+1,\za+1,\zb,\zi) + L_2(\zs,\za,\zb,\zi)
    \\ & \textstyle \null\qquad =  \frac{\zs-\za}{2\za-\zs+2} \sum_{j=0}^{\za} (-1)^{\za+j} \binom{j}{2\za-\zs+1}\binom{\za+1}{j}\binom{\zs-\za-\zb+j}{j-\zi}
        \\ & \textstyle \null\qquad\qquad + \frac{\zs-\za}{2\za-\zs+2} \binom{\za+1}{\zs-\za}\binom{\za-\zb+1}{\zi-\zb} - \binom{\za+1}{\zs-\za-1}\binom{\zs-\zb+1}{\zs-\za-\zb+\zi},
  \end{align*}
  which is zero by~\eqref{eq29835798}. Finally, suppose that \m{i=\beta}. Then 
    \begin{align*}
  \polyH_1'(-\beta) &\tempheightaer \textstyle = \smash{- \sum_{l=0}^\alpha \binom{l}{\beta}\Big(\binom{s-\alpha-\beta}{s-\alpha}^\partial \binom{l-\beta-1}{l+\alpha-s}+\binom{s-\alpha-\beta}{s-\alpha} \binom{l-\beta-1}{l+\alpha-s}^\partial\Big), } \\
    \polyH_2(-\beta) &\tempheightaer \textstyle = (-1)^{\beta+1}\binom{\beta-1}{s-\alpha-1}, \\
      \polyH_3'(-\beta) &\tempheightaer \textstyle =
    \smash{\sum_{j=0}^\alpha   (-1)^{\alpha+\beta+j+1} \binom{j}{2\alpha-s+1}\binom{\alpha+1}{j+1} \left(\binom{s-\alpha-\beta+j}{j-\beta}^\partial -\binom{s-\alpha-\beta+j}{s-\alpha}h_\beta\right)}, \\
      \polyH_4'(-\beta) &\tempheightaer =\textstyle  (-1)^\beta \binom{\alpha+1}{s-\alpha} (h_{\alpha-\beta} - h_\beta),
  \end{align*}
  where \m{h_t = 1 + \cdots + \frac1t}  is the harmonic number for \m{t \in \ZZ_{> 0}} and \m{h_{t}= 0} for \m{t \in \ZZ_{\LEQSLANT 0}}. 
  Since
  \begin{align*} & \textstyle \sum_{j=0}^\alpha   (-1)^{\alpha+\beta+j+1} \binom{j}{2\alpha-s+1}\binom{\alpha+1}{j+1} \binom{\polyr+s-\alpha+j}{s-\alpha} \\ & \textstyle \null\qquad =  (-1)^{\alpha+\beta} \left( \binom{\polyr+s-\alpha-1}{s-\alpha} - \binom{s-2\alpha-2}{s-\alpha}\right),\end{align*} we can simplify \m{ \polyH_3'(-\beta)} to
  \[\textstyle  \polyH_3'(-\beta) =  \sum_{j=0}^\alpha   (-1)^{\alpha} \binom{j}{2\alpha-s+1}\binom{\alpha+1}{j+1} \binom{\alpha-s-1}{j-\beta}^\partial -(-1)^{\beta}\left(\binom{\beta}{s-\alpha}-\binom{\alpha+1}{s-\alpha}\right)h_\beta.\]
  The equation \m{\polyH(-\beta)=0} is therefore equivalent to
  \[\textstyle L_3(s,\alpha,\beta) = R_3 (s,\alpha,\beta)\]
  with
  \begin{align*}
      \textstyle L_3 & \textstyle \defeq  \beta\polyH_1'(-\beta), \\
         R_3  & \textstyle \defeq (s-\alpha-\beta) (\polyH_3'(-\beta) + \polyH_4'(-\beta)) - \polyH_2(-\beta).
  \end{align*}
Let us show that \m{L_3(\zs,\za,\zb) = R_3(\zs,\za,\zb)} for all \m{\zs>\za\GEQSLANT \zb> 0}. For \m{\za=\zb} this is
  \[\textstyle (-1)^{\zs}\zb \left( \binom{\zb-1}{\zs-\zb}^\partial \binom{-1}{2\zb-\zs}-\binom{\zb-1}{\zs-\zb} \binom{-1}{2\zb-\zs}^\partial\right) =  (2\zb-\zs) \binom{\zb}{\zs-\zb} h_\zb + \binom{\zb-1}{\zs-\zb-1}. \]
  If \m{\zs>2\zb} then both sides are zero, if \m{\zs= 2\zb} then both sides are \m{1}, and if \m{2\zb > \zs > \zb} then both  sides are \m{\zb(h_{\zb-1}+\frac{1}{2\zb-\zs})}. Thus all we need to do is show that
  \[\textstyle R_3(\zs+1,\za+1,\zb)-R_3(\zs,\za,\zb)-L_3(\zs+1,\za+1,\zb)+L_3(\zs,\za,\zb) = 0\]
  for all \m{\zs>\za\GEQSLANT \zb> 0}. By using the equation
  \[\textstyle \sum_{j} (-1)^j \binom{j}{2\za-\zs+1} \binom{\za+1}{j} \binom{\polyr+\zs-\za-\zb+j}{j-\zb} = (-1)^{\za+1}  \binom{\za+1}{\zs-\za} \binom{\polyr + \za-\zb+1}{\za-\zb+1} \]
  we can get rid of the sum \m{\sum_{j}} and, after some simple algebraic manipulations, simplify this to
  \begin{align*}
      \textstyle \binom{\za+1}{\zb} \left( \binom{\zs-\za-\zb}{\zs-\za}^\partial \binom{\za-\zb}{2\za-\zs+1} + \binom{\zs-\za-\zb}{\zs-\za} \binom{\za-\zb}{2\za-\zs+1}^\partial \right)
     =  \frac{(-1)^{\zb+1}(\zs-\za-\zb)}{\zb(\za-\zb+1)} \binom{\za+1}{\zs-\za}. 
  \end{align*}
  We omit the full tedious details and just mention that since we are able to get rid of the sums \m{\sum_l} and \m{\sum_j} the aforementioned algebraic manipulations amount to simple cancellations.
  If \m{\zs\GEQSLANT \za+\zb} then 
   \begin{align*}
    &  \textstyle \binom{\za+1}{\zb} \left( \binom{\zs-\za-\zb}{\zs-\za}^\partial \binom{\za-\zb}{2\za-\zs+1} + \binom{\zs-\za-\zb}{\zs-\za} \binom{\za-\zb}{2\za-\zs+1}^\partial \right)
         \\ & \tempheightaer \textstyle \null\qquad = \frac{(-1)^{\zb+1} (\za+1)! (\za-\zb)! (\zs-\za-\zb)! (\zb-1)!}{\zb! (\za-\zb+1)! (2\za-\zs+1)! (\zs-\za-\zb-1)! (\zs-\za)!}
     \\ & \tempheightaer \textstyle \null\qquad =
      \frac{(-1)^{\zb+1}(\zs-\za-\zb) (\za+1)!}{\zb  (\za-\zb+1) (\zs-\za)! (2\za-\zs+1)!}
        =
      \frac{(-1)^{\zb+1}(\zs-\za-\zb)}{\zb(\za-\zb+1)} \binom{\za+1}{\zs-\za},
  \end{align*}
  and if \m{\zs < \za+\zb} then
     \begin{align*}
    &  \textstyle \binom{\za+1}{\zb} \left( \binom{\zs-\za-\zb}{\zs-\za}^\partial \binom{\za-\zb}{2\za-\zs+1} + \binom{\zs-\za-\zb}{\zs-\za} \binom{\za-\zb}{2\za-\zs+1}^\partial \right)
     \\ & \tempheightaer \textstyle \null\qquad = \frac{(-1)^{\zb} (\za+1)! (\zb-1)! (\za+\zb-\zs)! (\za-\zb)!}{\zb! (\za-\zb+1)! (\zs-\za)! (\za+\zb-\zs-1)! (2\za-\zs+1)!}
      \\ & \tempheightaer \textstyle \null\qquad = 
      \frac{(-1)^{\zb+1}(\zs-\za-\zb) (\za+1)!}{\zb  (\za-\zb+1) (\zs-\za)! (2\za-\zs+1)!}
        =
      \frac{(-1)^{\zb+1}(\zs-\za-\zb)}{\zb(\za-\zb+1)} \binom{\za+1}{\zs-\za}.
  \end{align*}
  We  have finally shown that if \m{i \in \{\beta,\ldots,\alpha\}} then \begin{align*}%
\textstyle \rowR_i (C_0,\ldots,C_\alpha)^T = 0.\end{align*}
\end{itemize}
  \end{proof}

  \section{Computing  \texorpdfstring{\m{\overline\Theta_{k,a}}}{\textTheta}}
    \label{secComp4}
    
    Throughout the proof we use the results from section~9 of~\cite{Ars1}, which we reproduce here without proofs for convenience.
    
    \begin{lemma}\label{LEMMAX2}
{
Suppose that \m{\alpha \in \{0,\ldots,\TEMPORARYI-1\}}.
\begin{enumerate}
    \item \label{111}
We have
  \begin{align*}
      \textstyle & \textstyle (T-a) \left(1 \sqbrackets{}{KZ}{\QQp}{ \theta^{n}x^{\alpha-n}y^{r-np-\alpha}}\right)
      \\ & \textstyle \null\qquad 
      = \sum_j (-1)^j \binom{n}{j} p^{j(p-1)+\alpha} \begin{Lmatrix}
      1 & 0 \\ 0 & p
      \end{Lmatrix} \sqbrackets{}{KZ}{\QQp}{x^{j(p-1)+\alpha}y^{r-j(p-1)-\alpha}} 
      \\ & \textstyle \null\qquad\qquad - a \sum_j (-1)^j \binom{n-\alpha}{j} \sqbrackets{}{KZ}{\QQp}{\theta^\alpha x^{j(p-1)}y^{r-j(p-1)-\alpha(p+1)}} + \BIGO{p^n}.
  \end{align*}
\item \label{222} The submodule \m{\Img(T-a)\subset\ind_{KZ}^G \TILDE\SIGMA_{r}} contains
  \begin{align*}
      & \textstyle \sum_{i} \left(\sum_{l=\beta-\gamma}^{\beta} C_l \binom{r-\beta+l}{i(p-1)+l}\right)\sqbrackets{}{KZ}{\QQp}{x^{i(p-1)+\beta}y^{r-i(p-1)-\beta}} \\ & \textstyle \null\qquad\qquad\qquad\qquad\qquad\qquad\qquad\qquad\qquad+ \BIGO{ap^{-\beta+v_C}+p^{p-1}}
  \end{align*}
  for all \m{0 \LEQSLANT \beta \LEQSLANT \gamma < \TEMPORARYI} and all families \m{\{C_l\}_{l\in \ZZ}} of elements of \m{\ZZ_p}, where
  \[\textstyle v_C = \min_{\beta-\gamma\LEQSLANT l \LEQSLANT \beta}(\PADICVALUATION(C_l)+l).\]
  The \m{\BIGO{ap^{-\beta+v_C}+p^{p-1}}} term is equal to  \m{\BIGO{p^{p-1}}} plus
  \[\textstyle -\frac{ap^{-\beta}}{p-1} \sum_{l=\beta-\gamma}^{\beta} C_{l}p^l\sum_{0\ne\MU\in\FF_p} [\MU]^{-l}\begin{Lmatrix}
p & [\MU] \\ 0 & 1
\end{Lmatrix}  \sqbrackets{}{KZ}{\QQp}{\theta^n  x^{\beta-l-n}y^{r-np-\beta+l}}.\]
\end{enumerate}
}
\end{lemma}

\begin{lemma}\label{LEMMAX3}
{
Suppose that \m{\alpha \in \ZZ} and \m{v \in \QQ} are such that
\begin{align*}\textstyle \alpha & \textstyle \in\{0,\ldots,\TEMPORARYI-1\}, \\ \textstyle v & \textstyle \LEQSLANT \PADICVALUATION(\TEMPORARYD_{\alpha}(D_\DARKCIRC)),\\v' \defeq \min\{\PADICVALUATION(a)-\alpha,v\} & \textstyle  \LEQSLANT \PADICVALUATION(\TEMPORARYD_{w}(D_\DARKCIRC)) \text{ for } \alpha < w < 2\TEMPORARYI-\alpha,
\\ v' & \textstyle <\PADICVALUATION(\TEMPORARYD_{w}(D_\DARKCIRC)) \text{ for } 0\LEQSLANT w<\alpha.
\end{align*} If, for \m{j \in \ZZ},
\[\textstyle\DELTAA_j \defeq (-1)^{j-1} (1-p)^{-\alpha} \binom{\alpha}{j-1} \TEMPORARYD_{\alpha}(D_\DARKCIRC),\]
then \m{v \LEQSLANT\PADICVALUATION(\TEMPORARYD_{\alpha}(\DELTAA_\DARKCIRC))\LEQSLANT \PADICVALUATION(\DELTAA_j)} for all \m{j \in \ZZ}, and
\begin{align*}
      \textstyle & \textstyle \sum_i (\DELTAA_i - D_i) \sqbrackets{}{KZ}{\QQp}{x^{i(p-1)+\alpha}y^{r-i(p-1)-\alpha}}
      \\ & \textstyle \null\qquad 
      =  \VERIFYIF{\alpha\LEQSLANT s}(-1)^{n+1} D_{\frac{r-s}{p-1}}
 \sqbrackets{}{KZ}{\QQp}{\theta^n x^{r-np-s+\alpha}y^{s-\alpha-n}} 
            \\ & \textstyle \null\qquad \qquad - D_0
 \sqbrackets{}{KZ}{\QQp}{\theta^n x^{\alpha-n}y^{r-np-\alpha}}  %
  \\ & \textstyle \null\qquad \qquad \qquad + E \sqbrackets{}{KZ}{\QQp}{\theta^{\alpha+1}h} + %
  F \sqbrackets{}{KZ}{\QQp}{h'} + \mathrm{ERR}_1 + \mathrm{ERR}_2,
  \end{align*}
  for some \m{\mathrm{ERR}_1} and \m{\mathrm{ERR}_2} such that
  \[\textstyle \mathrm{ERR}_1 \in\Img(T-a) \text{ and }\mathrm{ERR}_2 = \BIGO{p^{\TEMPORARYI-\PADICVALUATION(a)+v}+p^{\TEMPORARYI-\alpha}},\]
 some polynomials %
 \m{h} and \m{h'}, and some %
  \m{E,F\in\QQp} such that %
  \m{\PADICVALUATION(E)\GEQSLANT v'} and \m{\PADICVALUATION(F)>v'}.
}
\end{lemma}

\begin{lemma}\label{LEMMAX4}
{
Let  \m{\{C_l\}_{l\in \ZZ}} be any family of elements of \m{\ZZ_p}. Suppose that 
\m{\alpha \in \{0,\ldots, \TEMPORARYI-1\}} and \m{v \in \QQ}, and suppose that the constants
 \[\textstyle D_i \defeq %
 \VERIFYIF{i=0} C_{-1} + \VERIFYIF{0<i(p-1)<r-2\alpha}\sum_{l=0}^\alpha C_{l} \binom{r-\alpha+l}{i(p-1)+l}\]
satisfy the conditions of lemma~\ref{LEMMAX3}, \IE 
\begin{align*}\textstyle  \textstyle v & \textstyle \LEQSLANT \PADICVALUATION(\TEMPORARYD_{\alpha}(D_\DARKCIRC)),\\v' \defeq \min\{\PADICVALUATION(a)-\alpha,v\} & \textstyle  \LEQSLANT \PADICVALUATION(\TEMPORARYD_{w}(D_\DARKCIRC)) \text{ for } \alpha < w < 2\TEMPORARYI-\alpha,
\\ v' & \textstyle <\PADICVALUATION(\TEMPORARYD_{w}(D_\DARKCIRC)) \text{ for } 0\LEQSLANT w<\alpha.
\end{align*}
Moreover, suppose that \m{C_0} is a unit. Let
\[\textstyle\TEMPORARYD' \defeq (1-p)^{-\alpha}\TEMPORARYD_{\alpha}(D_\DARKCIRC) - C_{-1}   .\]
Suppose that \m{\PADICVALUATION(C_{-1})\GEQSLANT \PADICVALUATION(\TEMPORARYD')}.
\begin{enumerate}
    \item \label{7111} If \m{\PADICVALUATION(\TEMPORARYD') \LEQSLANT v'}
then there is some element  \m{\mathrm{gen}_1 \in \SCR I_a} that
represents a generator of \m{\FILT N_\alpha}. \\ \item \label{7222} If \m{\PADICVALUATION(a)-\alpha < v} then there is some element \m{\mathrm{gen}_2 \in \SCR I_a}  that represents a generator of a finite-codimensional submodule of 
\[\textstyle T\left(\ind_{KZ}^G\QUOTIENTI(\alpha)\right) =T\left(\FILT N_\alpha/ \ind_{KZ}^G\SUBMODULE(\alpha)\right),\]
 where \m{T} denotes the  endomorphism  of \m{\ind_{KZ}^G\QUOTIENTI(\alpha)} corresponding to the double coset of \m{\begin{Smatrix}
p & 0 \\ 0 & 1
\end{Smatrix}}.
\end{enumerate}
}
\end{lemma}

Let us now prove the following additional results.

    \begin{lemma}\label{LEMMATHM}
Let \m{\{C_l\}_{l\in \ZZ}} be any family of elements of~\m{\ZZ_p}. Suppose that \m{\alpha \in \ZZ} and \m{v \in \QQ} and the constants
 \[\textstyle D_i \defeq \VERIFYIF{i=0} C_{-1} + \VERIFYIF{0<i(p-1)<r-2\alpha}\sum_{l=0}^\alpha C_{l} \binom{r-\alpha+l}{i(p-1)+l}\] are such that
\begin{align*}\textstyle \alpha & \textstyle \in\{0,\ldots,\TEMPORARYI-1\}, \\ \textstyle v & \textstyle \LEQSLANT \PADICVALUATION(\TEMPORARYD_{\alpha}(D_\DARKCIRC)),
\\ v' \defeq \min\{\PADICVALUATION(a)-\alpha,v\}  & \textstyle <\PADICVALUATION(\TEMPORARYD_{\TEMPORARYAAAA}(D_\DARKCIRC)) \text{ for } 0\LEQSLANT \TEMPORARYAAAA<\alpha.
\end{align*} 
Let
\[\textstyle\TEMPORARYD' \defeq (1-p)^{-\alpha} \TEMPORARYD_{\alpha}(D_\DARKCIRC) - C_{-1}.\]
Then \m{\Img(T-a)} contains
\begin{align}& \textstyle (\TEMPORARYD'+C_{-1}) \sqbrackets{}{KZ}{\QQp}{\theta^\alpha x^{p-1}y^{r-\alpha(p+1)-p+1}}  + C_{-1}  \sqbrackets{}{KZ}{\QQp}{\theta^n x^{\alpha-n}y^{r-np-\alpha}}\nonumber
\\ & \textstyle\qquad \label{EQIMPX} + \sum_{\TEMPORARYA=\alpha+1}^{2\TEMPORARYI-\alpha-1} E_\TEMPORARYA 
 \sqbrackets{}{KZ}{\QQp}{\theta^{\TEMPORARYA}h_\TEMPORARYA %
 }  +  F \sqbrackets{}{KZ}{\QQp}{h'}
+ H%
,\end{align}
for some \m{h_\TEMPORARYA%
,h',E_\TEMPORARYA
,F,H} such that 
\begin{enumerate}
    \item\label{mainlemp1} \m{E_\TEMPORARYA = \TEMPORARYD_{\TEMPORARYA}(D_\DARKCIRC) + \BIGO{p^v} \cup \BIGO{\TEMPORARYD_{\alpha+1}(D_\DARKCIRC)} \cup\cdots\cup \BIGO{\TEMPORARYD_{\TEMPORARYA-1}(D_\DARKCIRC)}}, \\
     \item\label{mainlemp1a} if \m{\LESS{\TEMPORARYA+\alpha-s}\LEQSLANT\LESS{2\TEMPORARYA-s}\ne 0} then the reduction modulo \m{\maxideal} of \m{\textstyle \theta^{\TEMPORARYA}h_\TEMPORARYA} generates \m{N_\TEMPORARYA},\\
    \item\label{mainlemp2} \m{\PADICVALUATION(F)>v'}, and\\
    \item\label{mainlemp3} \m{H = \BIGO{p^{\TEMPORARYI-\PADICVALUATION(a)+v}+p^{\TEMPORARYI-\alpha}}} and if \m{\PADICVALUATION(a) -\alpha< v } then \[\textstyle \frac{1-p}{ap^{-\alpha}} H = g  \sqbrackets{}{KZ}{\QQp}{\theta^n x^{\alpha-n}y^{r-np-\alpha}} + \BIGO{p^{\TEMPORARYI-\PADICVALUATION(a)}}\]
    with
    \[\textstyle g = \sum_{\lambda\in\FF_p} C_0 \begin{Lmatrix}
    p & [\lambda] \\ 0 & 1\end{Lmatrix} + A \begin{Lmatrix}
    p & 0 \\ 0 & 1\end{Lmatrix} + \VERIFYIF{r\equiv_{p-1} 2\alpha} B  \begin{Lmatrix}
    0 & 1 \\ p & 0 \end{Lmatrix},\]
    where \[\textstyle A = -C_{-1}+\sum_{l=1}^\alpha C_l \binom{r-\alpha+l}{l}\] and \[\textstyle B=\sum_{l=0}^\alpha C_l \binom{r-\alpha+l}{s-\alpha}.\]
\end{enumerate}
\end{lemma}

\begin{proof}{Proof. } This lemma is essentially shown under a stronger hypothesis as lemma~\ref{LEMMAX4}. The stronger hypothesis consists of the three extra conditions that \m{\PADICVALUATION(\TEMPORARYD_{\TEMPORARYAAAA}(D_\DARKCIRC))\GEQSLANT \min\{\PADICVALUATION(a)-\alpha,v\}} for all \m{\alpha<\TEMPORARYAAAA<2\TEMPORARYI-\alpha}, that
\m{C_0\in \ZZ_p^\times}, and that \m{\PADICVALUATION(C_{-1})\GEQSLANT \PADICVALUATION(\TEMPORARYD')}.
These extra conditions are not used in the actual construction of the element in~\eqref{EQIMPX}, rather they are there to ensure that
\m{\PADICVALUATION(E_\TEMPORARYA) \GEQSLANT \min\{\PADICVALUATION(a)-\alpha,v\}} for all \m{\alpha<\TEMPORARYA<2\TEMPORARYI-\alpha}, that  the coefficient of \m{\begin{Smatrix}p & [\lambda] \\ 0 & 1\end{Smatrix}} in \m{g} is invertible, and that we get an integral element once we divide the element
\[\textstyle (\TEMPORARYD'+C_{-1}) \sqbrackets{}{KZ}{\QQp}{\theta^\alpha x^{p-1}y^{r-\alpha(p+1)-p+1}}  + C_{-1}  \sqbrackets{}{KZ}{\QQp}{\theta^n x^{\alpha-n}y^{r-np-\alpha}}\]
by \m{\TEMPORARYD'}.
Therefore we still get the existence of the element in~\eqref{EQIMPX} without these extra conditions, and to complete the proof of lemma~\ref{LEMMATHM} we need to verify the properties of \m{h_\TEMPORARYA,E_\TEMPORARYA,F,H,A}, and \m{B} claimed in~(\ref{mainlemp1}), (\ref{mainlemp1a}), (\ref{mainlemp2}), and (\ref{mainlemp3}).
 The \m{h_\TEMPORARYA} and \m{E_\TEMPORARYA} come from the proof of lemma~\ref{LEMMAX3}, and \m{E_\TEMPORARYA 
 \sqbrackets{}{KZ}{\QQp}{\theta^{\TEMPORARYA}h_\TEMPORARYA}} is
 \[\textstyle X_\TEMPORARYA 
 \sqbrackets{}{KZ}{\QQp}{\sum_{\lambda\ne 0} [-\lambda]^{r-\alpha-\TEMPORARYA} \begin{Lmatrix}1 & [\lambda] \\ 0 & 1\end{Lmatrix} (-\theta)^{n} x^{r-np-\TEMPORARYA}y^{\TEMPORARYA-n}},\]
with the notation for \m{X_\TEMPORARYA} from the proof of lemma~\ref{LEMMAX3}. Let \m{E_\TEMPORARYA=(-1)^{\TEMPORARYA+1}X_\TEMPORARYA}. Then condition~(\ref{mainlemp1}) is satisfied directly from the definition of \m{X_\TEMPORARYA}. Let
 \[\textstyle h_\TEMPORARYA = (-1)^{\TEMPORARYA+1} \sum_{\lambda\ne 0} [-\lambda]^{r-\alpha-\TEMPORARYA} \begin{Lmatrix}1 & [\lambda] \\ 0 & 1\end{Lmatrix} (-\theta)^n x^{r-np-\TEMPORARYA}y^{\TEMPORARYA-n}.\] This reduces modulo \m{\maxideal} to the element
\[\textstyle (-1)^{\TEMPORARYA} \sum_{\lambda\ne 0} [-\lambda]^{r-\alpha-\TEMPORARYA} \begin{Lmatrix}1 & [\lambda] \\ 0 & 1\end{Lmatrix} {Y^{\LESS{2\TEMPORARYA-r}}} = (-1)^{s-\alpha+1} \binom{\LESS{2\TEMPORARYA-s}}{\LESS{\TEMPORARYA+\alpha-s}} X^{\LESS{\TEMPORARYA+\alpha-s}} Y^{\LESS{\TEMPORARYA-\alpha}} \]
of  \[\textstyle \sigma_{\LESS{2\TEMPORARYA-r}}(r-\TEMPORARYA) \cong I_{r-2\TEMPORARYA}(\TEMPORARYA) / \sigma_{\MORE{r-2\TEMPORARYA}}(\TEMPORARYA)=\QUOTIENTI(\TEMPORARYA).\] This element is non-trivial and generates \m{N_\TEMPORARYA} if \m{\LESS{\TEMPORARYA+\alpha-s} \LEQSLANT \LESS{2\TEMPORARYA-s} \ne 0}, since then %
 \m{X^{\LESS{\TEMPORARYA+\alpha-s}} Y^{\LESS{\TEMPORARYA-\alpha}}} generates \m{N_\TEMPORARYA}. This verifies condition~(\ref{mainlemp1a}). Condition~(\ref{mainlemp2}) follows from the assumption \m{v'<\PADICVALUATION(\TEMPORARYD_{w}(D_\DARKCIRC))} for \m{0\LEQSLANT w < \alpha}, as in the proof of lemma~\ref{LEMMAX3}. Finally, condition~(\ref{mainlemp3}) follows from the description of the error term in lemma~\ref{LEMMAX2}, as in the proof of lemma~\ref{LEMMAX4}.
\end{proof}

\begin{corollary}\label{CORO1}
Let \m{\{C_l\}_{l\in \ZZ}} be any family of elements of~\m{\ZZ_p}. Suppose that \m{\alpha \in \{0,\ldots,\TEMPORARYI-1\}} and \m{v \in \QQ}, and suppose that the constants
 \[\textstyle D_i \defeq %
 \VERIFYIF{i=0} C_{-1} + \VERIFYIF{0<i(p-1)<r-2\alpha}\sum_{l=0}^\alpha C_{l} \binom{r-\alpha+l}{i(p-1)+l}\] are such that
\begin{align*}\textstyle  \textstyle v & \textstyle \LEQSLANT \PADICVALUATION(\TEMPORARYD_{\alpha}(D_\DARKCIRC)),
\\
v' \defeq \min\{\PADICVALUATION(a)-\alpha,v\}  & \textstyle \LEQSLANT \PADICVALUATION(\TEMPORARYD_{\TEMPORARYAAAA}(D_\DARKCIRC)) \text{ for } \alpha<\TEMPORARYAAAA<2\TEMPORARYI-\alpha,
\\ v'  & \textstyle <\PADICVALUATION(\TEMPORARYD_{\TEMPORARYAAAA}(D_\DARKCIRC)) \text{ for } 0\LEQSLANT w<\alpha.
\end{align*} 
Suppose also that \m{\PADICVALUATION(a)\not\in \ZZ}. Let \begin{align*} \textstyle\TEMPORARYD' &\textstyle\defeq (1-p)^{-\alpha} \TEMPORARYD_{\alpha}(D_\DARKCIRC) - C_{-1}, \\ \textstyle\veeC & \textstyle\defeq-C_{-1}+\sum_{l=1}^\alpha C_l \binom{r-\alpha+l}{l}.\end{align*}
 If \m{\star} then \m{\ast} is trivial modulo \m{\SCR I_a}, for each of the following pairs \[\textstyle(\star,\ast) = (\text{condition}, \text{representation}).\]

\begin{enumerate}
    \item\label{CORO1.1} \m{\Big(\PADICVALUATION(\TEMPORARYD') \LEQSLANT \min\{\PADICVALUATION(C_{-1}),v'\},\,\FILT N_\alpha\Big)}.\NEWLINEA
        \item\label{CORO1.2}  \m{\Big(v=\PADICVALUATION(C_{-1}) < \min\{\PADICVALUATION(\TEMPORARYD'),\PADICVALUATION(a)-\alpha\},\,\ind_{KZ}^G \SUBMODULE(\alpha)\Big)}.\NEWLINEA
            \item\label{CORO1.3}  \m{\Big(\PADICVALUATION(a)-\alpha < v\LEQSLANT \PADICVALUATION(C_{-1}) \annd \textstyle \veeC\in \ZZ_p^\times \annd C_0 \not\in\ZZ_p^\times \annd \LESS{2\alpha-r}> 0,\, \FILT N_\alpha\Big)}.\NEWLINEA
                                \item\label{CORO1.4} \m{\Big(\PADICVALUATION(a)-\alpha < v\LEQSLANT \PADICVALUATION(C_{-1})  \annd \textstyle %
   \veeC \in \ZZ_p^\times,\,\ind_{KZ}^G \QUOTIENTI(\alpha)\Big)}.\NEWLINEA
                \item\label{CORO1.5} \m{\Big(\PADICVALUATION(a)-\alpha < v\LEQSLANT \PADICVALUATION(C_{-1}) \annd C_0 \in \ZZ_p^\times,\,\mathbf{r}_1\Big)}, where \[\textstyle\mathbf{r}_1\] is a finite-codimensional submodule of \[\textstyle T(\ind_{KZ}^G \QUOTIENTI(\alpha)).\]
\end{enumerate}
\end{corollary}

\begin{proof}{Proof. } %
  There is one extra condition imposed in addition to the conditions from lemma~\ref{LEMMATHM}: that
  \[\textstyle v' \defeq \min\{\PADICVALUATION(a)-\alpha,v\}  \textstyle \LEQSLANT \PADICVALUATION(\TEMPORARYD_{\TEMPORARYAAAA}(D_\DARKCIRC)) \text{ for } \alpha<\TEMPORARYAAAA<2\TEMPORARYI-\alpha,\] and it ensures that
\m{\PADICVALUATION(E_\TEMPORARYA) \GEQSLANT v'} for all \m{\alpha<\TEMPORARYA<2\TEMPORARYI-\alpha}. Lemma~\ref{LEMMATHM} implies that the element in~\eqref{EQIMPX} is in \m{\Img(T-a)}. Let us call this element \m{\gamma}.

(\ref{CORO1.1}) 
The condition \m{\PADICVALUATION(\TEMPORARYD') \LEQSLANT \min\{\PADICVALUATION(C_{-1}),v'\}} ensures that if we divide \m{\gamma} by \m{\TEMPORARYD'} then the resulting element reduces modulo~\m{\maxideal} to a representative of a generator of \m{\FILT N_\alpha}.

(\ref{CORO1.2}) 
The condition \m{v=\PADICVALUATION(C_{-1}) < \min\{\PADICVALUATION(\TEMPORARYD'),\PADICVALUATION(a)-\alpha\}} ensures that if we divide \m{\gamma} by \m{C_{-1}} then the resulting element   reduces modulo~\m{\maxideal} to a representative of a generator of \m{\ind_{KZ}^G \SUBMODULE(\alpha)}.

(\ref{CORO1.3}, \ref{CORO1.4}, \ref{CORO1.5}) 
The condition \m{\PADICVALUATION(a)-\alpha < v\LEQSLANT \PADICVALUATION(C_{-1}) } ensures that the term with the dominant valuation in~\eqref{EQIMPX} is \m{H}, so we can divide \m{\gamma} by \m{ap^{-\alpha}} and obtain the element \m{L + \BIGO{p^{\TEMPORARYI-\PADICVALUATION(a)}}}, where \m{L} is defined by
\[\textstyle L \defeq  \left(\sum_{\LAMBDA\in\FF_p}C_0\begin{Lmatrix}
p & [\LAMBDA] \\ 0 & 1
\end{Lmatrix} + A \begin{Lmatrix}
p & 0 \\0 & 1
\end{Lmatrix} +\VERIFYIF{r\equiv_{p-1}2\alpha} B \begin{Lmatrix}
0 & 1 \\ p& 0
\end{Lmatrix}\right) \sqbrackets{}{KZ}{\QQp}{\theta^n x^{\alpha-n}y^{r-np-\alpha}}\] with \m{A} and \m{B} as  in lemma~\ref{LEMMATHM}. This element \m{L} is in \m{\Img(T-a)}, and it
 reduces modulo~\m{\maxideal} to a representative of
\[\textstyle \left(\sum_{\LAMBDA\in\FF_p}C_0\begin{Lmatrix}
p & [\LAMBDA] \\ 0 & 1
\end{Lmatrix} + A \begin{Lmatrix}
p & 0 \\0 & 1
\end{Lmatrix} +\VERIFYIF{r\equiv_{p-1}2\alpha} (-1)^{r-\alpha} B \begin{Lmatrix}
1 & 0 \\0 & p
\end{Lmatrix}\right)\sqbrackets{}{KZ}{\FFp}{X^{\LESS{2\alpha-r}}}.\] %
 As shown in the proof of lemma~\ref{LEMMAX4}, if \m{C_0\in \ZZ_p^\times} then this element always generates a finite-codimensional submodule of \[\textstyle T(\ind_{KZ}^G \QUOTIENTI(\alpha)),\] and if additionally 
 \m{A\ne 0} (over \m{\FF_p}) then in fact we have the stronger conclusion that it generates \[\textstyle \ind_{KZ}^G \QUOTIENTI(\alpha).\] Suppose on the other hand that \m{C_0=\BIGO{p}} and \m{A\in\ZZ_p^\times}. In that case we assume that \m{\LESS{2\alpha-r}>0} and therefore the reduction modulo \m{\maxideal} of \m{L} represents a generator of \m{\FILT N_\alpha}.\end{proof} %

\section{Proof of theorem~\ref{MAINTHM1a}}
\label{proof3}

 We prove theorem~\ref{MAINTHM1a} by proving nine propositions which give just enough information to conclude that \m{\overline \Theta_{k,a}} is irreducible, but not enough to classify it fully.

We assume that \[\textstyle r = s + \beta(p-1) + u_0p^t + \BIGO{p^{t+1}}\] for some \m{\beta \in \{0,\ldots,p-1\}}  and \m{u_0 \in \ZZ_p^\times} and \m{t \in \ZZ_{>0}}, and we write \m{\TEMPORARYJ = u_0p^t}.
As the main result of~\cite{Ars1} implies theorem~\ref{MAINTHM1a}
 for \m{s\GEQSLANT 2\TEMPORARYI}, we may assume that \[\textstyle s \in \{2,\ldots,2\TEMPORARYI-2\}.\]
Recall also that we assume
\m{\TEMPORARYI-1<\PADICVALUATION(a)<\TEMPORARYI} for some  \m{\TEMPORARYI\in\{1,\ldots,\frac{p-1}{2}\}},  and that \m{k>p^{100}} (and consequently \m{r > p^{99}}).

We now give a list of nine propositions, and show that their union implies theorem~\ref{MAINTHM1a}.

\begin{proposition}\label{propl1}
If \m{\alpha < \frac s2} then
\[\textstyle  %
 \left\{\begin{array}{rl} %
 \FILT N_\alpha  & \text{if }   \beta \in \{0,\ldots,\alpha-1\} \text{ and } \alpha > \PADICVALUATION(a)-t , \abcNL \ind_{KZ}^G\SUBMODULE(\alpha) & \text{otherwise}  \end{array}\right.\] is trivial modulo \m{\SCR I_a}.
\end{proposition}

\begin{proposition}\label{propl2}
If \m{\frac s2\LEQSLANT \alpha < s} and \m{\beta \not\in\{1,\ldots,\alpha+1\}} then \[\textstyle \FILT N_\alpha\] is trivial modulo \m{\SCR I_a}.
\end{proposition}

\begin{proposition}\label{propl3}
If \m{0 < \alpha < \frac s2} %
then \[\textstyle  %
 \left\{\begin{array}{rl} %
T(\ind_{KZ}^G\QUOTIENTI(\alpha)) & \text{if }  \beta \in %
 \{0,\ldots, \alpha\}  \text{ and } \alpha >\PADICVALUATION(a)-t, \abcNL \FILT N_{s-\alpha} %
&  \text{if } \beta \in %
 \{0,\ldots, \alpha\}   \text{ and }  \alpha <\PADICVALUATION(a)-t , \abcNL
\FILT N_\alpha & \text{if } \beta \in \{\alpha+1,\ldots,s-\alpha\}, \abcNL \FILT N_{s-\alpha} %
& \text{if } \beta > s-\alpha
 \end{array}\right.\]
is trivial modulo \m{\SCR I_a}.
\end{proposition}

\begin{proposition}\label{propl4}
If \m{\frac s2\LEQSLANT \alpha < s} and \m{(\alpha,\beta) \ne (\frac s2,\frac s2+1)} then
\[\textstyle  %
 \left\{\begin{array}{rl} %
T(\ind_{KZ}^G\QUOTIENTI(\alpha)) & \text{if }   \beta \in \{1,\ldots,s-\alpha\}  \text{ and } s-\alpha > \PADICVALUATION(a)-t, \abcNL
T(\ind_{KZ}^G\QUOTIENTI(\alpha)) & \text{if }   \beta \in \{s-\alpha+1,\ldots,\alpha\} \text{ and } \alpha > \PADICVALUATION(a)-t , \abcNL\FILT N_{\alpha} %
&  \text{otherwise} \end{array}\right.\]
is trivial modulo \m{\SCR I_a}.
\end{proposition}

\begin{proposition}\label{propl5}
If %
\m{\alpha \GEQSLANT s} then \[\textstyle  %
 \left\{\begin{array}{rl} %
T(\ind_{KZ}^G\QUOTIENTI(\alpha)) & \text{if } \alpha = \max\{\TEMPORARYI-t-1,\beta-1\}, \abcNL \FILT N_\alpha %
& \text{otherwise} \end{array}\right.\] is trivial modulo \m{\SCR I_a}. 
\end{proposition}

\begin{proposition}\label{propl6}
If \m{\beta \in \{1,\ldots,\frac s2-1\}} and \m{t > \TEMPORARYI-\frac s2-2} then  \[\textstyle \FILT N_{s/2+1}\] %
 is trivial modulo \m{\SCR I_a}. \end{proposition}

\begin{proposition}\label{propl7}
If \m{\beta \in \{1,\ldots,\frac{s}{2}-1\}} and 
\m{t = \TEMPORARYI-\frac{s}{2}} then \[\textstyle \FILT N_{s/2-1}\] %
 is trivial modulo \m{\SCR I_a}.
\end{proposition}

\begin{proposition}\label{propl8}
If \m{\beta \in \{\frac{s}{2},\frac{s}{2}+1\}} and 
\m{t > \TEMPORARYI-\frac s2-1} then \[\textstyle \FILT N_{s/2+1}\] %
 is trivial modulo \m{\SCR I_a}.
\end{proposition}

\begin{proposition}\label{propl9}
If \m{\beta =\frac{s}{2}+1} and 
\m{t = \TEMPORARYI-\frac{s}{2}-1} then \[\textstyle \ind_{KZ}^G\SUBMODULE(\frac s2+1)\] %
 is trivial modulo \m{\SCR I_a}.
\end{proposition}

\begin{proof}{Proof that propositions~\ref{propl1}--\ref{propl9} imply theorem~\ref{MAINTHM1a}. }
Let us assume that \m{\overline{\Theta}_{k,a}} is reducible with the goal of reaching a contradiction. The classification given by theorem 2 in~\cite{Ars1} implies that \m{\overline{\Theta}_{k,a}} has two infinite-dimensional factors, each of which is a quotient of a representation in the set
 \[\{\ind_{KZ}^G\SUBMODULE(\alpha) \,|\, 0\LEQSLANT\alpha<\TEMPORARYI\} \cup \{\ind_{KZ}^G\QUOTIENTI(\alpha) \,|\, 0\LEQSLANT\alpha<\TEMPORARYI\},\]
 and moreover that the following classification is true.
 \begin{enumerate}
     \item If the two representations are \m{\ind_{KZ}^G\SUBMODULE(\alpha_1)} and \m{\ind_{KZ}^G\SUBMODULE(\alpha_2)} then \[\textstyle\alpha_1+\alpha_2\equiv_{p-1} s+1.\]
     \item If the two representations are \m{\ind_{KZ}^G\SUBMODULE(\alpha_1)} and \m{\ind_{KZ}^G\QUOTIENTI(\alpha_2)} then \[\textstyle\alpha_1-\alpha_2\equiv_{p-1} 1.\]
     \item If the two representations are \m{\ind_{KZ}^G\QUOTIENTI(\alpha_1)} and \m{\ind_{KZ}^G\QUOTIENTI(\alpha_2)} then \[\textstyle\alpha_1+\alpha_2\equiv_{p-1} s-1.\]
 \end{enumerate}
 The facts that \begin{align*}\alpha_1+\alpha_2 & \textstyle\in \{0,\ldots,2\TEMPORARYI-2\} \subseteq \{0,\ldots,p-3\}, \\\alpha_1-\alpha_2 & \textstyle\in \{1-\TEMPORARYI,\ldots,\TEMPORARYI-1\} \subseteq \{-\frac{p-3}{2},\ldots,\frac{p-3}{2}\}, \\
 s &\in \{2,\ldots,2\TEMPORARYI-2\} \subseteq \{2,\ldots,p-3\}\end{align*} imply that the following classification is true as well.
 \begin{enumerate}
     \item If the two representations are \m{\ind_{KZ}^G\SUBMODULE(\alpha_1)} and \m{\ind_{KZ}^G\SUBMODULE(\alpha_2)} then \[\textstyle\alpha_1+\alpha_2 = s+1.\]
     \item If the two representations are \m{\ind_{KZ}^G\SUBMODULE(\alpha_1)} and \m{\ind_{KZ}^G\QUOTIENTI(\alpha_2)} then \[\textstyle\alpha_1=\alpha_2+1.\]
     \item If the two representations are \m{\ind_{KZ}^G\QUOTIENTI(\alpha_1)} and \m{\ind_{KZ}^G\QUOTIENTI(\alpha_2)} then \[\textstyle\alpha_1+\alpha_2=s-1.\]
 \end{enumerate}
 This classification and %
 propositions~\ref{propl1}, \ref{propl2}, \ref{propl3}, \ref{propl4}, and~\ref{propl5} together imply that one of the two representations must be either \m{\ind_{KZ}^G \SUBMODULE(\frac s2)} or \m{\ind_{KZ}^G \QUOTIENTI(\frac s2)}, and in that case the other representation is either \[\textstyle\textstyle\ind_{KZ}^G \SUBMODULE(\frac s2+1)\] (which can only happen if \m{\beta \in \{1,\ldots,\frac s2-1\}} and \m{t > \TEMPORARYI-\frac s2} or \m{\beta\in\{\frac s2,\frac s2+1\}} and \m{t > \TEMPORARYI-\frac s2-2}), or 
\[\textstyle \textstyle \ind_{KZ}^G \QUOTIENTI(\frac s2-1)\] (which can only happen if \m{s=2} or \m{\beta \in \{1,\ldots,\frac s2-1\}} and \m{t = \TEMPORARYI-\frac s2}).
 In the latter case if \m{s=2} then either \m{1\sqbrackets{}{KZ}{\QQp}{x^2y^{r-2}} \in \SCR I_a} generates \m{\ind_{KZ}^G\QUOTIENTI(0)}, or \m{\TEMPORARYI\LEQSLANT 2} in which case \m{\overline{V}_{k,a}} is known to be irreducible.
Propositions~\ref{propl4}, \ref{propl6}, \ref{propl7}, \ref{propl8}, and~\ref{propl9}
 exclude all of the remaining possibilities. Thus if we assume that  \m{\overline{\Theta}_{k,a}} is reducible we reach a contradiction, so \m{\overline{\Theta}_{k,a}} must be irreducible.
\end{proof}

\emph{Proof of proposition~\ref{propl1}. }
First suppose that \m{\beta\GEQSLANT \alpha}. We apply part~(\ref{CORO1.2}) of corollary~\ref{CORO1} with \m{v =0} and
\begin{align*} 
 C_j = \left\{\begin{array}{rl}  (-1)^{\alpha}\binom{s-r}{\alpha} & \text{if } j = -1,\Cjnewline
 0 & \text{if } j = 0,\Cjnewline
 (-1)^{\alpha-j}\binom{s-\alpha+1}{\alpha-j} & \text{if } j \in \{1,\ldots,\alpha\}.
 \end{array}\right.\end{align*}
 Since \[\textstyle \binom{s-r}{\alpha} = \binom{\beta}{\alpha}+\BIGO{p} \in \ZZ_p^\times,\] the two conditions we need to verify %
  are \m{\PADICVALUATION(\TEMPORARYD_{w}(D_\DARKCIRC)) > 0} for \m{0\LEQSLANT w < \alpha} and \m{\PADICVALUATION(\TEMPORARYD')>0}. These two conditions are equivalent to the system of equations 
 \begin{align}\label{eqn4235}\nonumber& \textstyle  \textstyle \textstyle  \sum_{j=1}^\alpha (-1)^{\alpha-j} \binom{s-\alpha+1}{\alpha-j} \sum_{i>0} \binom{r-\alpha+j}{i(p-1)+j}\binom{i(p-1)}{w} \\ &  \textstyle\null\qquad\qquad\qquad\qquad\qquad\qquad = (-1)^\alpha(\VERIFYIF{w=\alpha}-\VERIFYIF{w=0}) \binom{s-r}{\alpha}+\BIGO{p} \end{align} for \m{0\LEQSLANT w\LEQSLANT \alpha}.
Let \m{F_{w,j}(\polyr,\polys) \in \FF_p[\polyr,\polys]} denote the polynomial defined in lemma~\ref{lemmaComb555}. Since
\[\textstyle \sum_{i>0} \binom{r-\alpha+j}{i(p-1)+j}\binom{i(p-1)}{w} = F_{w,j}(r,s)\]
by~\eqref{ceqn6}, the conclusion of that lemma when evaluated at \m{\polyr=r} and \m{\polys=s} implies~\eqref{eqn4235}.
 Thus if \m{\beta\GEQSLANT\alpha} then we can apply part~(\ref{CORO1.2})  of corollary~\ref{CORO1} and conclude that  \m{\textstyle \ind_{KZ}^G\SUBMODULE(\alpha)} is trivial modulo \m{\SCR I_a}.
 
 Suppose now that \m{\beta \in \{0,\ldots,\alpha-1\}}.
If \m{t>\PADICVALUATION(a)-\alpha} then the proof of theorem 17 in~\cite{Ars1} applies here nearly verbatim %
 since \[\textstyle\binom{s-\alpha+1}{\alpha} \in \ZZ_p^\times,\] and in fact we can conclude that \m{\FILT N_\alpha} is trivial modulo \m{\SCR I_a}. %
  So let us suppose that \m{t<\PADICVALUATION(a)-\alpha}. %
   We apply part~(\ref{CORO1.2})  of corollary~\ref{CORO1} with \m{v = t} and
 \begin{align*} 
 C_j = \left\{\begin{array}{rl}  (-1)^{\alpha}\binom{s-r}{\alpha} & \text{if } j = -1,\Cjnewline
 0 & \text{if } j = 0,\Cjnewline
 (-1)^{\alpha-j}\binom{s-\alpha+1}{\alpha-j} + pC_j^\ast & \text{if } j \in \{1,\ldots,\alpha\},
 \end{array}\right.\end{align*}
 for some constants \m{C_1^\ast,\ldots,C_\alpha^\ast} yet to be chosen. Clearly \[\textstyle \PADICVALUATION(C_{-1})=t < \PADICVALUATION(a)-\alpha,\] and the other conditions that need to be satisfied in order for corollary~\ref{CORO1} to be applicable are \begin{align*}t & \textstyle < \PADICVALUATION(\TEMPORARYD'), \\ t & \textstyle \LEQSLANT\PADICVALUATION(\TEMPORARYD_{w}(D_\DARKCIRC)) \text{ for } \alpha\LEQSLANT w<2\TEMPORARYI-\alpha, \\ t & \textstyle <\PADICVALUATION(\TEMPORARYD_{w}(D_\DARKCIRC)) \text{ for } 0\LEQSLANT w<\alpha.\end{align*} Let us consider the matrix \m{A = (A_{w,j})_{0\LEQSLANT w,j \LEQSLANT\alpha}} that has integer entries
\begin{align*}
     A_{w,j}  & \textstyle = \sum_{0<i(p-1)<r-2\alpha} \binom{r-\alpha+j}{i(p-1)+j}\binom{i(p-1)}{w}.
\end{align*}
 Then exactly as in the proof of theorem 17 in~\cite{Ars1} we can show that
\[\textstyle A = S + \epsilon N + \BIGO{\epsilon p},\] where
\begin{align*}
    S_{w,j} & = \textstyle \sum_{i=1}^\beta \binom{s+\beta(p-1)-\alpha+j}{i(p-1)+j}\binom{i(p-1)}{w},\\
    N_{w,j} & = \textstyle \sum_v (-1)^{w-v} \binom{j+w-v-1}{w-v} \binom{s+\beta(p-1)-\alpha+j}{v}^\partial \sum_{i=0}^\beta \binom{s+\beta(p-1)-\alpha+j-v}{i(p-1)+j-v}\\ & \textstyle\null\qquad -\VERIFYIF{w=0} \binom{s+\beta(p-1)-\alpha+j}{j}^\partial.
\end{align*}
We still have equation~\ref{eqn4235} since the constants are the same, and since \[\textstyle\binom{s-r}{\alpha}=\BIGO{p},\] we have
\[\textstyle S (C_0,\ldots,C_\alpha)^T = (\BIGO{p},\ldots,\BIGO{p})^T.\]
Let
\m{B = B_\alpha} be the 
\m{(\alpha+1)\times(\alpha+1)} matrix defined in lemma~\ref{lemmaComb2}. That lemma implies that 
\m{B} encodes precisely the row operations that transform \m{S} into a matrix with zeros outside the rows indexed \m{1,\ldots,\beta} and such that \[\textstyle (BS)_{w,j} = p^{-\VERIFYIF{j=0}}\binom{s+\beta(p-1)-\alpha+j}{w(p-1)+j}\]
when \m{w \in\{1,\ldots,\beta\}}.
 We thus have
\[\textstyle BS (C_0,\ldots,C_\alpha)^T = (0,\BIGO{p},\ldots, \BIGO{p},0,\ldots)^T,\]
where the only entries of the vector on the right that can possibly be non-zero are the ones indexed \m{1,\ldots,\beta}. As in the proof of theorem 17 in~\cite{Ars1} we note that \m{S} has rank \m{\beta} and therefore we can choose \m{C_1^\ast,\ldots,C_\alpha^\ast} in a way that \m{(C_0,\ldots,C_\alpha)^T\in \Ker BS}. Then  \m{\TEMPORARYD_{w}(D_\DARKCIRC) = \BIGO{\epsilon}} for all \m{w}, and the conditions that need to be satisfied are \m{\TEMPORARYD_{w}(D_\DARKCIRC) = \BIGO{\epsilon p}} for \m{0\LEQSLANT w<\alpha} and \m{\TEMPORARYD'=\BIGO{\epsilon p}}. %
 These two conditions are equivalent to the single equation
 \[\textstyle A(C_0,\ldots,C_\alpha)^T = (-C_{-1},0,\ldots,0,C_{-1}) + \BIGO{\epsilon p},\]
 which is itself equivalent to
 \begin{align*} & \textstyle BN (C_0,\ldots,C_\alpha)^T \\ & \null\qquad\textstyle = \left(0,-\binom{\alpha}{1}(C_{-1}\epsilon^{-1}),\ldots,(-1)^\alpha \binom{\alpha}{\alpha}(C_{-1}\epsilon^{-1})\right)^T + BSv + \BIGO{p}\end{align*}
 for some \m{v}. %
 Thus, if \m{\overline{R}} is the \m{\alpha\times\alpha} matrix over \m{\FF_p} obtained from \m{\overline{BN}} by replacing the rows indexed \m{1,\ldots,\beta} with the corresponding rows of \m{\overline{BS}} and then discarding the zeroth row and the zeroth column, the condition that needs to be satisfied is equivalent to the claim that
 \[\textstyle \left(-(1-\VERIFYIF{1\LEQSLANT\beta})\binom{\alpha}{1},\ldots,(-1)^\alpha (1-\VERIFYIF{\alpha\LEQSLANT\beta}) \binom{\alpha}{\alpha}\right)^T\]
is in the image of \m{\overline{R}} (since \m{C_0=\BIGO{p}} and \m{C_{-1}\epsilon^{-1} \in \ZZ_p^\times}). This is indeed the case since \m{\overline R} is the lower right  \m{\alpha\times\alpha} submatrix of the matrix \m{\overline Q} defined in the proof of theorem 17 in~\cite{Ars1} (where it is shown that \m{\overline{Q}} is equal to the matrix \m{M} from lemma~\ref{lemmaComb4}) and is therefore upper triangular with units on the diagonal. Thus we can apply part~(\ref{CORO1.2})  of corollary~\ref{CORO1} with \m{v = t} and conclude that \m{\ind_{KZ}^G\SUBMODULE(\alpha)} is trivial modulo \m{\SCR I_a}.%
\hspace{\stretch{1}}     {\BLACKSQUARE}

\emph{Proof of proposition~\ref{propl2}. }
Let us define \m{\polyconst_{-1}(\polyr),\ldots,\polyconst_\alpha(\polyr) \in \ZZ_p[\polyr]} as
\begin{align*} 
 \polyconst_j(\polyr) = \left\{\begin{array}{rl} %
\binom{s-\polyr-1}{\alpha+1} & \text{if } j=-1, \\[5pt] \binom{\alpha}{s-\alpha-1}^{-1}\frac{s-\polyr}{\alpha+1} %
& \text{if } j=0,
\\[5pt] \frac{(-1)^{j+1}}{j+1} \binom{s-\alpha-1}{\alpha-j} (\polyr-\alpha) & \text{if } j \in \{1,\ldots,\alpha\}. \end{array}\right.\end{align*}
We apply part~(\ref{CORO1.1})  of corollary~\ref{CORO1} with \m{v =0} and
\[\textstyle (C_{-1},C_0,\ldots,C_\alpha) = (\polyconst_{-1}(r),\polyconst_0(r),\ldots,\polyconst_{\alpha}(r)).\]
 The two conditions we need to verify %
  are \m{\PADICVALUATION(\TEMPORARYD_{w}(D_\DARKCIRC)) > 0} for \m{0\LEQSLANT w < \alpha} and \m{\PADICVALUATION(\TEMPORARYD')=0}. These two conditions follow from the system of equations \begin{align}\label{eq1234}\textstyle \sum_{j=0}^\alpha C_j \sum_{0 < i(p-1)<r-2\alpha}  \binom{r-\alpha+j}{i(p-1)+j}\binom{i(p-1)}{w} = -\VERIFYIF{w=0} \binom{s-r-1}{\alpha+1} + \BIGO{p}\end{align} for \m{0\LEQSLANT w\LEQSLANT \alpha}. 
Let \m{F_{w,j}(\polyr) \in \FF_p[\polyr]} denote the polynomial
\[\textstyle \sum_v (-1)^{w-v} \binom{j+w-v-1}{w-v}\binom{\polyr-\alpha+j}{v} \binom{\MORE{s-\alpha+j-v}}{\LESS{j-v}} - \binom{\polyr-\alpha+j}{j}\binom{0}{w} - \binom{\polyr-\alpha+j}{s-\alpha}\binom{\polyr-s}{w}.\]
By~\eqref{ceqn6}, \[\textstyle \sum_{0 < i(p-1)<r-2\alpha} \binom{r-\alpha+j}{i(p-1)+j}\binom{i(p-1)}{w} = F_{w,j}(r),\]
 so the conclusion of lemma~\ref{lemmaComb5} evaluated at~\m{\polyr=r} implies~\eqref{eq1234}. Thus
we can apply part~(\ref{CORO1.1}) of corollary~\ref{CORO1} and conclude that \m{\FILT N_\alpha} is trivial modulo \m{\SCR I_a}.
\hspace{\stretch{1}}     {\BLACKSQUARE}

\emph{Proof of proposition~\ref{propl3}. }
First let us assume that \m{\beta \in \{0,\ldots,\alpha\}}. If we attempt to copy the proof of theorem 17 in~\cite{Ars1} in this setting, the one place where we run into problems is that some entries of the extended associated matrix \m{N} are not integers (\IE when we extend the number of rows in \m{A}, \m{S}, and \m{N} to \m{2\TEMPORARYI-\alpha} by defining \m{A_{w,j}}, \m{S_{w,j}}, and \m{N_{w,j}} with the same equations used for the first \m{\alpha+1} rows, we get entries which are not integers). To be more specific, the equation for \m{N_{w,0}} in this setting is
 \[\textstyle pN_{w,0} = \binom{s+\beta(p-1)-\alpha}{w}^\partial \sum_{i>0} \binom{s+\beta(p-1)-\alpha-w}{i(p-1)-w} + \BIGO{p}, \]
 where the second term is \m{\BIGO{p}} because it is still true that \[\textstyle %
 \sum_{i>0} \binom{r-\alpha-w}{i(p-1)-w} - \sum_{i>0} \binom{s+\beta(p-1)-\alpha-w}{i(p-1)-w}%
  = \BIGO{\epsilon p}.\]
 On the other hand,
 \begin{align*}\textstyle \sum_{i>0} \binom{s+\beta(p-1)-\alpha-w}{i(p-1)-w} & \textstyle = \sum_{l=0}^w (-1)^l \binom{w}{l} \sum_{i>0} \binom{s+\beta(p-1)-\alpha-l}{i(p-1)}
 \\
 & \textstyle = (-1)^{s-\alpha} \binom{w}{s-\alpha} + \BIGO{p}.
 \end{align*}
 So \m{A_{w,0} = S_{w,0} + \BIGO{\epsilon}} is integral if \m{w < s-\alpha} and
 \[\textstyle A_{w,0} = S_{w,0} + (-1)^{s-\alpha} \binom{w}{s-\alpha}  \binom{s-\alpha-\beta}{w}^\partial \epsilon p^{-1} + \BIGO{\epsilon}\] if \m{w\GEQSLANT s-\alpha}.
Note that \m{\beta \in \{0,\ldots,\alpha\}} and \m{s > 2\alpha} by assumption, so \m{S_{w,0}} is still always integral, and if \m{s-\alpha \LEQSLANT w < 2\TEMPORARYI-\alpha} then
\[\textstyle \binom{s-\alpha-\beta}{w}^\partial =  \frac{(-1)^{s-\alpha-\beta-w+1}}{w \binom{w-1}{s-\alpha-\beta}} \in \ZZ_p^\times.\]
What this means is that if we proceed with the proof of theorem 17 in~\cite{Ars1} and apply lemma~\ref{LEMMATHM} with the  constants \m{(C_{-1},C_0,\ldots,C_\alpha)} constructed there such that \m{C_0} is a unit, then
we obtain  an element
\begin{align*}& \textstyle (\TEMPORARYD'+C_{-1}) \sqbrackets{}{KZ}{\QQp}{\theta^\alpha x^{p-1}y^{r-\alpha(p+1)-p+1}}  + C_{-1}  \sqbrackets{}{KZ}{\QQp}{\theta^n x^{\alpha-n}y^{r-np-\alpha}}\nonumber
\\ & \textstyle\qquad  + \sum_{\TEMPORARYA=\alpha+1}^{2\TEMPORARYI-\alpha-1} E_\TEMPORARYA 
 \sqbrackets{}{KZ}{\QQp}{\theta^{\TEMPORARYA}h_\TEMPORARYA %
 }  +  F \sqbrackets{}{KZ}{\QQp}{h'}
+ H\end{align*}
which is in \m{\Img(T-a)} and is such that \begin{align*}
    \textstyle \PADICVALUATION (C_{-1}) & \textstyle = \PADICVALUATION(\TEMPORARYD') = t+1, \\ \textstyle \PADICVALUATION(E_\TEMPORARYA) & \textstyle \GEQSLANT t+1 \text{ for } \alpha+1 \LEQSLANT \TEMPORARYA < s-\alpha,\\
    \PADICVALUATION(F) & \textstyle > t+1,
\end{align*} and with \m{H} as in lemma~\ref{LEMMATHM}. However, \m{\PADICVALUATION(E_{s-\alpha}) = t} and \m{\PADICVALUATION(E_\TEMPORARYA) \GEQSLANT t} for \m{\TEMPORARYA > s-\alpha}. Therefore if \m{t > \PADICVALUATION(a)-\alpha} then the dominant term is \m{H} and we can conclude that a submodule of finite codimension in \m{T(\ind_{KZ}^G\QUOTIENTI(\alpha)) } is trivial modulo \m{\SCR I_a}, and if \m{t < \PADICVALUATION(a) - \alpha} then the dominant term is \[\textstyle E_{s-\alpha}\sqbrackets{}{KZ}{\QQp}{\theta^{s-\alpha}h_{s-\alpha}} \] and hence \m{\FILT N_{s-\alpha}} is trivial modulo \m{\SCR I_a} by part~(\ref{mainlemp1a}) of lemma~\ref{LEMMATHM}.

Now let us assume that \m{\beta > \alpha}. We use the constants constructed in the second bullet point of the proof of theorem 17 in~\cite{Ars1}, and we apply lemma~\ref{LEMMATHM}.
This gives  an element
\begin{align*}& \textstyle \TEMPORARYD' \sqbrackets{}{KZ}{\QQp}{\theta^\alpha x^{p-1}y^{r-\alpha(p+1)-p+1}}  %
\\ & \textstyle\qquad  + \sum_{\TEMPORARYA=\alpha+1}^{2\TEMPORARYI-\alpha-1} E_\TEMPORARYA 
 \sqbrackets{}{KZ}{\QQp}{\theta^{\TEMPORARYA}h_\TEMPORARYA %
 }  +  F \sqbrackets{}{KZ}{\QQp}{h'}
+ H\end{align*}
which is in \m{\Img(T-a)} and is such that \begin{align*}
    \PADICVALUATION(\TEMPORARYD')&=1,\\
    \PADICVALUATION(E_\TEMPORARYA) &\GEQSLANT 1 \text{ for }\alpha+1 \LEQSLANT \TEMPORARYA < s-\alpha,\\
    \PADICVALUATION(F) &> 1,\\
    \textstyle \PADICVALUATION(E_{s-\alpha}) & \textstyle = \PADICVALUATION((r-\alpha)_{s-\alpha}), %
\end{align*} and with \m{H} as in lemma~\ref{LEMMATHM}. This time the dominant term is either \[\textstyle \TEMPORARYD' \sqbrackets{}{KZ}{\QQp}{\theta^\alpha x^{p-1}y^{r-\alpha(p+1)-p+1}}\] or \[\textstyle E_{s-\alpha}\sqbrackets{}{KZ}{\QQp}{\theta^{s-\alpha}h_{s-\alpha}} \]  depending on whether \m{\beta \in \{\alpha+1,\ldots,s-\alpha\}} or \m{\beta > s-\alpha}. Thus in the former case \m{\FILT N_{\alpha}} is trivial modulo \m{\SCR I_a}, and in the latter case \m{\FILT N_{s-\alpha}} is trivial modulo \m{\SCR I_a}.
\hspace{\stretch{1}}     {\BLACKSQUARE}

\emph{Proof of proposition~\ref{propl4}. }
By proposition~\ref{propl2} we may assume that \m{\beta\not\in \{1,\ldots,\alpha+1\}}, and by proposition~\ref{propl3} we may assume that \m{\beta \ne \alpha+1}. If \m{\alpha\ne \frac s2} and \m{\beta \in \{1,\ldots,s-\alpha\}} and  \m{s-\alpha < \PADICVALUATION(a)-t} then the claim follows from proposition~\ref{propl3}. Thus it is enough to show that if \m{\beta \in \{1,\ldots,\alpha\}} then \[\textstyle  %
 \left\{\begin{array}{rl} %
T(\ind_{KZ}^G\QUOTIENTI(\alpha)) & \text{if } \alpha > \PADICVALUATION(a)-t , \\\FILT N_{\alpha} %
&  \text{if } \alpha < \PADICVALUATION(a)-t  \end{array}\right.\]
is trivial modulo \m{\SCR I_a}. 
If \m{\alpha < \PADICVALUATION(a)-t} we apply part~(\ref{CORO1.1}) of corollary~\ref{CORO1}, and if \m{\alpha > \PADICVALUATION(a)-t}) we apply part~(\ref{CORO1.5}) of corollary~\ref{CORO1}. In both cases we choose  \m{v=t} and 
\begin{align*} 
 C_j = \left\{\begin{array}{rl} \frac{(-1)^{\alpha+\beta}(s-\alpha)(\alpha-\beta+1)}{\beta^2(2\alpha-s+1) \binom{\alpha}{\beta}} \binom{\alpha}{s-\alpha} \epsilon & \text{if } j=-1, \\[7pt] 1
& \text{if } j=0,
\\[5pt] \frac{(-1)^{j+1}(s-\alpha-\beta)}{\beta} \binom{j}{2\alpha-s+1} \binom{\alpha+1}{j+1} & \text{if } j \in \{1,\ldots,\alpha\}. \end{array}\right.\end{align*}
 Since \m{\PADICVALUATION(C_{-1}) = t} and \m{C_0 = 1}, the conditions we need to verify in order to be able to apply corollary~\ref{CORO1} are
 \begin{align*}
     t & \LEQSLANT\PADICVALUATION(\TEMPORARYD_{w}(D_\DARKCIRC)) \text{ for }\alpha\LEQSLANT w<2\TEMPORARYI-\alpha,\\
     t & <\PADICVALUATION(\TEMPORARYD_{w}(D_\DARKCIRC)) \text{ for } 0\LEQSLANT w<\alpha, \\
     \TEMPORARYD' & = -C_{-1} + \BIGO{\epsilon p}.
 \end{align*}
Let us consider the matrix \[\textstyle A = (A_{w,j})_{0\LEQSLANT w,j \LEQSLANT\alpha}\] that has integer entries
\begin{align*}
     A_{w,j}  & \textstyle = \sum_{0<i(p-1)<r-2\alpha} \binom{r-\alpha+j}{i(p-1)+j}\binom{i(p-1)}{w}.
\end{align*}
Then the second and third conditions are equivalent to
the claim that
\[\textstyle A(C_0,\ldots,C_\alpha)^T = \left(-C_1+\BIGO{\epsilon p},\BIGO{\epsilon p},\ldots,\BIGO{\epsilon p}\right)^T.\]
 As in the proof of the approximation claim in the proof of the main result of~\cite{Ars1} (and as in proposition~\ref{propl1}) we can show that
\[\textstyle A = S + \epsilon N + \BIGO{\epsilon p},\] where
\begin{align*}
    S_{w,j} & = \textstyle \sum_{i=1}^{\beta} \binom{s+\beta(p-1)-\alpha+j}{i(p-1)+j}\binom{i(p-1)}{w}-\binom{s+\beta(p-1)-\alpha+j}{s-\alpha}\binom{\beta(p-1)}{w},\\
    N_{w,j} & = \textstyle \sum_v (-1)^{w-v} \binom{j+w-v-1}{w-v} \binom{s+\beta(p-1)-\alpha+j}{v}^\partial \sum_{i=0}^\beta \binom{s+\beta(p-1)-\alpha+j-v}{i(p-1)+j-v}\\ & \textstyle\null\qquad -\VERIFYIF{w=0} \binom{s+\beta(p-1)-\alpha+j}{j}^\partial - \binom{s-\alpha-\beta+j}{s-\alpha}^\partial\binom{-\beta}{w}-\binom{s-\alpha-\beta+j}{s-\alpha}\binom{-\beta}{w}^\partial.
\end{align*} The first condition follows from an argument similar to the one in the fourth bullet point in the proof of theorem 17 in~\cite{Ars1}: if we extend the number of rows in \m{A}, \m{S}, and \m{N} to \m{2\TEMPORARYI-\alpha} by defining \m{A_{w,j}}, \m{S_{w,j}}, and \m{N_{w,j}} with the same equations used for the first \m{\alpha+1} rows, then we have \m{A \equiv S \bmod \epsilon} and so we can replace \m{A} with \m{S + \BIGO{\epsilon}}, and \m{\TEMPORARYD_{w}(D_\DARKCIRC)} for each \m{\alpha\LEQSLANT w<2\TEMPORARYI-\alpha} is a \m{\ZZ_p}-linear combination of \m{\TEMPORARYD_{0}(D_\DARKCIRC)=\BIGO{\epsilon},\ldots,\TEMPORARYD_{\alpha}(D_\DARKCIRC)=\BIGO{\epsilon}}.
 And, as in the proof of theorem 17 in~\cite{Ars1}, the second and third conditions follow if \begin{align*}
S (C_{0},\ldots,C_{\alpha})^T & = 0, \\
\textstyle N (C_{0},\ldots,C_{\alpha})^T  & \textstyle= \left(-C_1\epsilon^{-1},0,\ldots,0\right)^T + Sv + \BIGO{p}\text{ for some }v.\end{align*}  Let
\m{B = B_\alpha} be the 
\m{(\alpha+1)\times(\alpha+1)} matrix defined in lemma~\ref{lemmaComb2}. 
 Then \m{BS} has zeros outside of the rows indexed \m{1,\ldots,\beta-1}, and \[\textstyle (BS)_{i,j} = \binom{s+\beta(p-1)-\alpha+j}{i(p-1)+j}\]
for \m{i \in \{1,\ldots,\beta-1\}}. Let \m{\overline{R}} denote the \m{(\alpha+1)\times(\alpha+1)} matrix over \m{\FF_p} obtained from \m{\overline{BN}} by replacing the rows indexed \m{1,\ldots,\beta-1} with the corresponding rows of \m{\overline{BS}}. As in the proof of theorem 17 in~\cite{Ars1} we can compute 
\begin{align*} \textstyle (\overline{BN})_{i,j} & \tempheightaer\textstyle = \sum_{l,v=0}^\alpha (-1)^{i+l+v} \binom{l}{i} \binom{j-v}{l-v} \binom{s-\alpha-\beta+j}{v}^\partial \binom{\MORE{s-\alpha+j-v}}{\LESS{j-v}} %
 \\ & \tempheightaer \textstyle\null\qquad - \VERIFYIF{i=0} \binom{s-\alpha-\beta+j}{j}^\partial - \VERIFYIF{i=\beta} \binom{s-\alpha-\beta+j}{s-\alpha}^\partial \\ & \tempheightaer\textstyle\null\qquad
  - (-1)^{i}  \binom{s-\alpha-\beta+j}{s-\alpha} \sum_{l=0}^\alpha \binom{l}{i} \binom{l-\beta-1}{l}^\partial. \end{align*} Thus lemma~\ref{comblemma1} implies that
\[\textstyle \overline{R} (C_{0},C_1,\ldots,C_\alpha)^T = 
\left(\frac{(-1)^{\alpha+\beta+1}(s-\alpha)(\alpha-\beta+1)}{\beta^2(2\alpha-s+1) \binom{\alpha}{\beta}} \binom{\alpha}{s-\alpha},0,\ldots,0\right)^T.
\] So the conditions we need to apply corollary~\ref{CORO1} are indeed satisfied, and that completes the proof.
\hspace{\stretch{1}}     {\BLACKSQUARE}

\emph{Proof of proposition~\ref{propl5}. }
This is the first time that we consider an \m{\alpha} such that \m{\alpha\GEQSLANT s}. The major difference in this scenario is that \m{s} is not the ``correct'' remainder of \m{r} to work with and instead we should consider the number that is congruent to \m{r} mod \m{p-1} and belongs to the set \m{\alpha+1,\ldots,p-\alpha-1}. Let us therefore define \m{s_\alpha = \MORE{r-\alpha}+\alpha}, and in particular let us note that \m{s_\alpha = s} for \m{s>\alpha} (which has hitherto always been the case). 
Then the computations in the proof of theorem 17 in~\cite{Ars1} work out exactly the same if we replace every instance of \m{s} with \m{s_\alpha} (and the restricted sum ``\m{\sum_{i>0}}'' with ``\m{\sum_{0<i(p-1)<r-\alpha}}'' when \m{s_\alpha = p-1}). The sufficient condition for these computations to work is
 \[\textstyle \binom{s_\alpha-\alpha}{2\TEMPORARYI-\alpha} \in \ZZ_p^\times,\] which is indeed the case since
 \m{s_\alpha-\alpha = p-1+s-\alpha \GEQSLANT 2\TEMPORARYI-\alpha}.
 So there is an analogous version of theorem 17 in~\cite{Ars1}, and we can conclude the desired result---as the proof of theorem 17 in~\cite{Ars1} works nearly without modification, we omit the full details of the arguments.
\hspace{\stretch{1}}     {\BLACKSQUARE}

\emph{Proof of proposition~\ref{propl6}. }
Let us write \m{\alpha= \frac s2+1} and, as the claim we want to prove is vacuous for \m{s=2}, let us assume that \m{s\GEQSLANT 4} and in particular \m{\alpha\GEQSLANT3}. We  apply part~(\ref{CORO1.3})  of corollary~\ref{CORO1} with \m{v} chosen in the open interval \m{(\PADICVALUATION(a)-\alpha,t)} and
\begin{align*} 
 C_j = \left\{\begin{array}{rl}   0 & \text{if } j  \in \{-1,0\}, \Cjnewline
 (-1)^{j}\binom{\alpha-2}{j} + (-1)^{j+1} (\alpha-2) \binom{\alpha-2}{j-1} + pC_j^\ast & \text{if } j \in \{1,\ldots,\alpha\},
 \end{array}\right.\end{align*} for some constants \m{C_1^\ast,\ldots,C_\alpha^\ast} yet to be chosen. The conditions necessary for the lemma to be applicable are satisfied if 
 \m{\veeC = \textstyle \sum_{j} C_j \binom{r-\alpha+j}{j} \in \ZZ_p^\times} and 
 \[\TEMPORARYD_w(D_\DARKCIRC) = \BIGO{\epsilon}\]
 for \m{0\LEQSLANT w < 2\TEMPORARYI-\alpha}. We have
 \begin{align*}
     \textstyle \textstyle \veeC & \tempheightaer\textstyle = \sum_{j} C_j \binom{s-\alpha-\beta+j}{s-\alpha-\beta} + \BIGO{p} \\ & \tempheightaer\textstyle =
-1 + \sum_{j} \left((-1)^{j}\binom{\alpha-2}{j} + (-1)^{j+1} (\alpha-2) \binom{\alpha-2}{j-1}\right) \binom{s-\alpha-\beta+j}{s-\alpha-\beta} + \BIGO{p} \\ & \tempheightaer\textstyle = -1 + \BIGO{p} \in \ZZ_p^\times \end{align*}
by~\eqref{ceqn5}  since \m{\alpha-2>s-\alpha-\beta}. And, since  \[\textstyle j\LEQSLANT s-\alpha-\beta+j \LEQSLANT s-\beta < p-i,\] we also have \[\textstyle \binom{s+\beta(p-1)-\alpha+j}{i(p-1)+j} = \binom{\beta}{i} \binom{s-\alpha-\beta+j}{j-i} + \BIGO{p}.\] Thus the equality \m{\TEMPORARYD_w(D_\DARKCIRC) = \BIGO{p}} follows from the fact that
\[\textstyle \sum_{j} (-1)^j \binom{\alpha-2}{j-2} \binom{s-\alpha-\beta+j}{s-\alpha-\beta+i} = 0,\]
which follows from~\eqref{ceqn5} since \m{\textstyle \alpha-2>\alpha-2-\beta+i=s-\alpha-\beta+i}. Moreover, we can choose \[\textstyle C_1^\ast,\ldots,C_\alpha^\ast\] in a way that \m{\TEMPORARYD_w(D_\DARKCIRC)  = 0} for \m{0\LEQSLANT w < 2\TEMPORARYI - \alpha} similarly as in the proof of theorem 17 in~\cite{Ars1} since the reduction modulo~\m{p} of the matrix
     \[\textstyle \left(\binom{s+\beta(p-1)-\alpha+j}{i(p-1)+j}\right)_{1\LEQSLANT i,j< \beta} = \left(\binom{\beta}{i}\binom{s-\alpha-\beta+j}{j-i} + \BIGO{p}\right)_{1\LEQSLANT i,j< \beta}\]
     is upper triangular with units on the diagonal.
 Thus the conditions we need to apply corollary~\ref{CORO1} are satisfied and we can conclude that \m{\FILT N_{s/2+1}} is trivial modulo \m{\SCR I_a}.
\hspace{\stretch{1}}     {\BLACKSQUARE}

\emph{Proof of proposition~\ref{propl7}. }
Let us write \m{\alpha= \frac s2-1} %
 and, as the claim we want to prove is vacuous for \m{s=2}, let us assume that \m{s\GEQSLANT 4} and in particular \m{\alpha\GEQSLANT3}. The only obstruction in the proof of proposition~\ref{propl3} that prevents us from concluding that \m{\FILT N_{s/2-1}} is trivial modulo \m{\SCR I_a} is that the dominant terms are \[\textstyle \textstyle E_{\TEMPORARYA}\sqbrackets{}{KZ}{\QQp}{\theta^{\TEMPORARYA}h_{\TEMPORARYA}}\] for \m{\frac s2 < \TEMPORARYA \LEQSLANT 2\TEMPORARYI-\frac s2}
 rather than \m{H}. %
  We can see from proposition~\ref{propl6} that \[\textstyle E_{s/2+1}\sqbrackets{}{KZ}{\QQp}{\theta^{s/2+1}h_{s/2+1}} = x_1 + x_2,\] with \m{\PADICVALUATION(x_1)\GEQSLANT t+1}, and with \m{x_2 \in \Img(T-a)}. Since the valuation of the coefficient of \m{H} is less than \m{t+1}, we can remove the obstruction coming from \[\textstyle E_{s/2+1}\sqbrackets{}{KZ}{\QQp}{\theta^{s/2+1}h_{s/2+1}}\] by replacing it with \m{x_1}. If \m{s=2\TEMPORARYI-2} then this is the only obstruction and we can  conclude that \m{\FILT N_{s/2-1}} is trivial modulo \m{\SCR I_a}. Now suppose that \m{s < 2\TEMPORARYI-2}. Then just as in the proof of theorem 17 in~\cite{Ars1} we  can apply part~(\ref{CORO1.1}) of corollary~\ref{CORO1} and conclude that \m{\FILT N_\alpha} is trivial modulo \m{\SCR I_a} as long as \m{(\epsilon,0,\ldots,0)^T} is in the image of
  the matrix \m{A = (A_{w,j})_{0\LEQSLANT w,j \LEQSLANT\alpha}} that has integer entries
\begin{align*}
     A_{w,j}  & \textstyle = \sum_{i>0} \binom{r-\alpha+j}{i(p-1)+j}\binom{i(p-1)}{w} = S_{w,j} + \epsilon N_{w,j} + \BIGO{\epsilon p}
\end{align*}
with \m{S} and \m{N} as in proposition~\ref{propl1}.
However, this time we can deduce more than that: since \m{s < 2\TEMPORARYI-2} %
 it follows that
 \[\textstyle 1\sqbrackets{}{KZ}{\QQp}{x^{s/2+1}y^{r-s/2-1}}\] is equal to \[\textstyle g_1 \sqbrackets{}{KZ}{\QQp}{\theta^{s/2+1} x^{s/2-n+1}y^{r-np-s/2-1}} + x_3\] for some \m{g_1} with \m{\PADICVALUATION(g_1)\GEQSLANT \PADICVALUATION(a)-\frac s2-1} and some \m{x_3 \in\Img(T-a)}. This in turn by proposition~\ref{propl6} is equal to \[\textstyle g_2 \sqbrackets{}{KZ}{\QQp}{\theta^{s/2+1}h_2} + x_4\] for some \m{g_2} with \m{\PADICVALUATION(g_2)\GEQSLANT t}, some \m{h_2}, and some \m{x_4\in \Img(T-a)}. Here we use the fact that the valuation of the constant \m{C_1} from proposition~\ref{propl6} is at least one and therefore the corresponding term \m{H} is \[\textstyle g_2 \sqbrackets{}{KZ}{\QQp}{\theta^{s/2+1} x^{s/2-n+1}y^{r-np-s/2-1}}+ x_5 + \BIGO{\epsilon} \]  for some \m{g_3} with \m{\PADICVALUATION(g_3)= \PADICVALUATION(a)-\frac s2-1} and some \m{x_5\in \Img(T-a)}. In general the error term would be \[\textstyle C_1ap^{-s/2} g_4\sqbrackets{}{KZ}{\QQp}{\theta^{s/2} x^{s/2-n}y^{r-np-s/2}}+\BIGO{\epsilon}\] rather than \m{\BIGO{\epsilon}}---a description of this error term is given in part~(\ref{222}) of lemma~\ref{LEMMAX2}. 
 This implies that we can add a constant multiple of \[\textstyle 1\sqbrackets{}{KZ}{\QQp}{x^{s/2+1}y^{r-s/2-1}}\] to the element \[\textstyle \sum_i D_i \sqbrackets{}{KZ}{\QQp}{x^{i(p-1)+\alpha}y^{r-i(p-1)-\alpha}} + \BIGO{ap^{-\alpha}}\] from the proof of lemma~\ref{LEMMAX4}, and we can translate this back to adding the extra column \[\textstyle \left(\binom{r-\alpha}{s-\alpha},\ldots,\binom{r}{s}\right)^T\] to \m{A}.
As in proposition~\ref{propl1} we can then reduce showing that  \m{(\epsilon,0,\ldots,0)^T} is in the image of \m{A} to showing that \[\textstyle (1,0,\ldots,0)^T\] is in the image of the \m{(\alpha+1)\times(\alpha+2)} matrix \m{\overline{R}} which is obtained from the matrix \m{\overline{Q}} defined in the proof of theorem 17 in~\cite{Ars1} by replacing all entries in the first row with zeros (because this time we do not divide the corresponding row of \m{A} by \m{p}) and by adding an extra column corresponding to the extra column of \m{A}. Thus, if we index the extra column to be the zeroth column, the lower right \m{\alpha\times\alpha} submatrix of \m{\overline R} is upper triangular with units on the diagonal, the first column of \m{\overline{R}} is identically zero, and all entries of the first row of \m{\overline{R}} except for \m{\overline R_{0,0}} are zero. As when computing \m{(\overline{BN})_{i,j}} in proposition~\ref{propl4} we can find that
\[\textstyle \overline{R}_{0,0} = \sum_{l=0}^\alpha %
 \binom{l-\beta-1}{l}^\partial = \polyH'(-\beta-1)\]
 with
 \[\textstyle \polyH(\polyr) = \sum_{l=0}^\alpha \binom{\polyr+l}{l} = \binom{\polyr+\alpha+1}{\alpha}.\]
 Thus 
 \[\textstyle  \overline{R}_{0,0} = \binom{\alpha-\beta}{\alpha}^\partial = \frac{(-1)^{\beta+1}}{\beta \binom{\alpha}{\beta}} \ne 0,\]
which implies that \m{(1,0,\ldots,0)^T} is in the image of \m{\overline{R}}. %
 Thus the conditions we need to apply corollary~\ref{CORO1} are satisfied and we can conclude that \m{\FILT N_{s/2-1}} is trivial modulo \m{\SCR I_a}.
\hspace{\stretch{1}}     {\BLACKSQUARE}

\emph{Proof of proposition~\ref{propl8}. }
Let us write \m{\alpha= \frac s2+1}. 
The reason why the proof of proposition~\ref{propl6} does not work for \m{\beta \in \{\alpha-1,\alpha\}} is because \m{\veeC=\BIGO{p}} for the constructed constants \m{C_j}. 
However, since \m{t>\PADICVALUATION(a)-\frac s2}, if \m{\veeC \in p\ZZ_p^\times} then the dominant term coming from lemma~\ref{LEMMATHM} is
\[\textstyle H = b_H \begin{Lmatrix} p & 0 \\ 0 & 1 \end{Lmatrix} \sqbrackets{}{KZ}{\QQp}{\theta^n x^{\alpha-n}y^{r-np-\alpha}} + \BIGO{p^{\TEMPORARYI-\alpha+1}}\] for the constant \[\textstyle b_H = \frac{ap^{-\alpha}}{1-p}\veeC\] which has valuation  \m{\PADICVALUATION(a)-\alpha+1}. As in proposition~\ref{propl7} it is crucial here that \m{C_1 = \BIGO{p}}. Just as in the proof of  proposition~\ref{propl6} we can reduce the claim we want to show to proving that there exist constants \m{C_1,\ldots,C_\alpha \in \ZZ_p} such that \m{C_1 = \BIGO{p}} and
 \[\textstyle \left(\binom{s+\beta(p-1)-\alpha+j}{i(p-1)+j}\right)_{0\LEQSLANT i< \beta, 0<j\LEQSLANT\alpha} (C_1,\ldots,C_\alpha)^T = (p,0,\ldots,0)^T.\]
Therefore it is enough to show that the square matrix
\[\textstyle A_0 =  \left(p^{\VERIFYIF{j=1}-\VERIFYIF{i\LEQSLANT\beta-\alpha+1}}\binom{s+\beta(p-1)-\alpha+j}{i(p-1)+j}\right)_{0\LEQSLANT i< \beta,\, \alpha-\beta<j\LEQSLANT\alpha}\]
has integer entries and is invertible (over \m{\ZZ_p}), as then we can recover 
\begin{align*} 
 \left\{\begin{array}{rl}  C_1 = 0 \text{ and } (C_2,\ldots,C_\alpha)^T = A_0^{-1} (1,0,\ldots,0)^T & \text{if } \beta=\alpha-1,\\
 (C_1/p,\ldots,C_\alpha)^T = A_0^{-1} (1,0,\ldots,0)^T & \text{if } \beta=\alpha.
 \end{array}\right.\end{align*}
 This follows from lemma~\ref{comblemma1cc}.
 So the conditions we need to apply corollary~\ref{CORO1} are satisfied and we can conclude that \m{\FILT N_{s/2+1}} is trivial modulo \m{\SCR I_a}.
\hspace{\stretch{1}}     {\BLACKSQUARE}

\emph{Proof of proposition~\ref{propl9}. }
Let us write \m{\alpha= \frac s2+1}. This time the proofs of both parts~(\ref{propl6}) and~(\ref{propl8}) break down since \m{\veeC = \BIGO{p}} %
 and the dominant term is no longer \m{H}. Let us slightly tweak these constants and instead use
\begin{align*} 
 C_j = \left\{\begin{array}{rl}   (-1)^\alpha \epsilon & \text{if } j =-1, \Cjnewline (-1)^{\alpha+j+1}\alpha\binom{\alpha-2}{j-2} & \text{if } j  \in \{0,\ldots,\alpha\}.
 \end{array}\right.\end{align*}
  Let \m{\overline{R}} be the matrix constructed in proposition~\ref{propl4}. Then just as in the proof of proposition~\ref{propl6} we can show that \m{\veeC = \BIGO{p}}, and just as in the proof of proposition~\ref{propl1} we can show that the dominant term coming from equation~\eqref{EQIMPX} in lemma~\ref{LEMMATHM} is
\[\textstyle (\TEMPORARYD'+C_{-1}) \sqbrackets{}{KZ}{\QQp}{\theta^\alpha x^{p-1}y^{r-\alpha(p+1)-p+1}}  + C_{-1}  \sqbrackets{}{KZ}{\QQp}{\theta^n x^{\alpha-n}y^{r-np-\alpha}}\] (and therefore that \m{\ind_{KZ}^G\SUBMODULE(\frac s2+1)} %
 is trivial modulo \m{\SCR I_a}) as long as
 \begin{align*}
    \overline{R} (C_0,\ldots,C_\alpha)^T = (0,\ldots,0,1)^T.
 \end{align*}
 This follows from lemma~\ref{comblemma1p}.
     Thus the conditions we need to apply corollary~\ref{CORO1} are satisfied and we can conclude that \m{\ind_{KZ}^G \SUBMODULE(\frac s2+1)} is trivial modulo \m{\SCR I_a}.
\hspace{\stretch{1}}     {\BLACKSQUARE}

\end{document}